\documentclass{tac}
\usepackage[all]{xy}
\usepackage{amsmath,amssymb}
\title{Butterflies in a semi-abelian context}
\author{O. Abbad, S. Mantovani, G. Metere, E. M. Vitale}
\address
{Dipartimento di matematica\\
Universit\`a degli studi di Milano\\
Via C. Saldini 50\\
20133 Milano, Italia ;\\
Institut de recherche en math\'ematique et physique\\
Universit\'e catholique de Louvain\\
Chemin du Cyclotron 2\\
B 1348 Louvain-la-Neuve, Belgique}
\copyrightyear{2010}
\thanks{Financial support: FNRS grant 1.5.276.09 is
gratefully acknowledged.}
\eaddress{oabbad@hotmail.com\\
sandra.mantovani@unimi.it\\
giuseppe.metere@unimi.it\\
enrico.vitale@uclouvain.be}

\newcommand\ev{\text{\rm {ev}}}
\newcommand\coker{\text{\rm {coker\,}}}
\newcommand\Aut{\text{\rm {Aut}}}
\newcommand\Ext{\text{\it {Ext}}}
\newcommand\Der{\text{\it {Der}}}
\newcommand\Rng{\text{\it {Rng}}}
\newcommand\Grp{\text{\it {Grp}}}
\newcommand\Lie{\text{\it {Lie}}}
\newcommand\Set{\text{\it {Set}}}
\newcommand\Vect{\text{\it {Vect}}}

\newcommand\Hom{\text{\it {Hom}}}
\newcommand\GrpdC{\text{\it {Grpd$(\cC)$}}}
\newcommand\CatC{\text{\it {Cat$(\cC)$}}}
\newcommand\XModC{\text{\it {XMod$(\cC)$}}}
\newcommand\Grpd{\text{\it {Grpd}}}
\newcommand{\A}{{\mathbb A}}
\newcommand{\B}{{\mathbb B}}
\newcommand{\C}{{\mathbb C}}
\newcommand{\D}{{\mathbb D}}
\newcommand{\E}{{\mathbb E}}
\newcommand{\G}{{\mathbb G}}
\newcommand{\bH}{{\mathbb H}}
\newcommand{\bK}{{\mathbb K}}

\newcommand{\X}{{\mathbb X}}
\newcommand{\K}{{\mathbb K}}
\newcommand{\cA}{{\mathcal A}}
\newcommand{\cB}{{\mathcal B}}
\newcommand{\cC}{{\mathcal C}}

\newcommand{\cF}{{\mathcal F}}
\newcommand{\cI}{{\mathcal I}}
\newcommand{\cP}{{\mathcal P}}

\newtheorem{Theorem}{Theorem}[section]
\newtheorem{Lemma}[Theorem]{Lemma}
\newtheorem{Proposition}[Theorem]{Proposition}
\newtheorem{Definition}[Theorem]{Definition}
\newtheorem{Corollary}[Theorem]{Corollary}
\newtheorem{Remark}[Theorem]{Remark}

\newtheorem{Notation}[Theorem]{Notation}

\newtheorem{IntGroupoid}[Theorem]{The 2-category of internal groupoids}
\newtheorem{DiscFib}[Theorem]{Discrete (co)fibrations}
\newtheorem{WeakEquiv}[Theorem]{Weak equivalences}
\newtheorem{XMod}[Theorem]{Internal crossed modules}
\newtheorem{Farfalle}[Theorem]{Internal butterflies}
\newtheorem{Torsori}[Theorem]{Fractors}
\newtheorem{CompFarfalle}[Theorem]{Composition of butterflies}
\newtheorem{IdFarfalle}[Theorem]{Identity butterflies}
\newtheorem{SplitFarfalle}[Theorem]{From morphisms to split butterflies}
\newtheorem{SplitFarfalleToMor}[Theorem]{From split butterflies to morphisms}
\newtheorem{CompRid}[Theorem]{Reduced composition}
\newtheorem{HomF}[Theorem]{The universal homomorphism}
\newtheorem{Classification}[Theorem]{Classification}
\newtheorem{FreeEx}[Theorem]{Free exact categories}
\newtheorem{Replacement}[Theorem]{Replacement}
\newtheorem{BicatFract}[Theorem]{Bicategories of fractions}

\newtheorem{SemiAbCat}[Theorem]{Protomodular and semi-abelian categories}
\newtheorem{nBarrKock}[Theorem]{(normalized) Barr-Kock Theorem}
\newtheorem{HuqIsSmith}[Theorem]{The ``Huq = Smith'' condition}
\newtheorem{InternalAct}[Theorem]{Internal object actions}
\newtheorem{SpanButterfly}[Theorem]{Span associated to a Butterfly}
\newtheorem{2CatXMods}[Theorem]{The 2-category of crossed modules}
\newtheorem{CoopXMods}[Theorem]{The morphisms $\kappa$ and $\iota$ cooperate}

\begin{document}

\maketitle

\begin{abstract}
It is known that monoidal functors between internal groupoids in the category $\Grp$ of groups
constitute the bicategory of fractions of the 2-category $\Grpd(\Grp)$ of internal groupoids,
internal functors and internal natural transformations in $\Grp$ with respect to weak equivalences.
Monoidal functors can be described equivalently by a kind of weak morphisms
introduced by B. Noohi under the name of ``butterflies''.
In order to internalize monoidal functors in a wide context, we
introduce the notion of  internal butterflies between internal crossed modules in a semi-abelian
category $\cC$, and we show that they are morphisms of a bicategory $\cB(\cC).$
Our main result states that, when in $\cC$ the notions of Huq commutator and Smith
commutator  coincide,
then the bicategory $\cB(\cC)$ of internal butterflies is the bicategory of fractions
of $\Grpd(\cC)$ with respect to weak equivalences (that is, internal functors which
are internally fully faithful and essentially surjective on objects).
\end{abstract}

\section{Introduction}

A  groupoid in the category of groups is a special case of strict monoidal
category, tensor product being provided by the group structure on objects and
arrows. Therefore, beyond internal functors, as arrows between  groupoids
in groups we can consider monoidal functors, that is functors between the underlying
categories
$$F \colon \bH \to \G$$
equipped with a natural and coherent family of isomorphisms
$$F^{x,y} \colon Fx + Fy \to F(x+y) \;\;\;\;\; x,y \in H_0$$
Both notions of monoidal functor and internal functor are relevant as morphisms
of  groupoids in groups (just to cite an example, as special case of
monoidal functors we get group extensions, whereas in the same case internal
functors give split extensions, see Section 6), so the question of expressing in
an internal way monoidal functors arises.

Three progresses have been recently accomplished in this direction.
In \cite{Noohi05} (see also \cite{AN09}) B. Noohi has proved that the
bicategory having  groupoids in groups as objects and monoidal functors as 1-cells
can be equivalently
described using crossed modules of groups as objects and what he calls
``butterflies'' as arrows. Moreover, in a paper with E. Aldrovandi \cite{AN09}, the theory is
pushed forward
in order to include the more general situation where groups are replaced by internal groups in a
(Grothendieck) topos. Noohi's butterflies (of \cite{Noohi05}) has been studied in the case of
Lie algebras on
a field in \cite{Abbad10}, and in \cite{NoohiLie}, where it is proved that butterflies
between differential
crossed modules (i.e.\ crossed modules of Lie algebras) represent homomorphisms
of strict Lie 2-algebras.

On the other hand, in \cite{Vitale10} it has been proved that the
bicategory of  groupoids in groups and monoidal functors is the bicategory of fractions of
the bicategory
of  groupoids in groups and internal functors with respect to weak equivalences.
Once again, the same result holds replacing groups with Lie algebras and monoidal
functors with homomorphisms of strict Lie 2-algebras.
In \cite{Dupont08}, M. Dupont has proved that butterflies provides the bicategory of fractions
of internal functors
with respect to weak equivalences when working internally to any abelian category.

The aim of this paper is to unify the results in \cite{Noohi05}, \cite{Abbad10}, \cite{NoohiLie},
\cite{Dupont08} and \cite{Vitale10}.
We introduce and study the bicategory $\mathcal{B}(\mathcal{C})$ of crossed modules
and  butterflies in a semi-abelian category $\mathcal{C}$.
The main result  is Theorem \ref{ThMain}, where we prove that $\mathcal{B}(\mathcal{C})$
is the bicategory of
fractions of the bicategory of  groupoids and internal functors with respect to weak equivalences.
This result gives a general answer to the specific problem recalled above:
it describes   weak internal functors that generalize at once monoidal functors in \Grpd(\Grp) and
homomorphisms in \Grpd(\Lie), and works for other 2-dimensional algebraic settings, as for
groupoids of Leibniz algebras, associative algebras, rings etc.
The non pointed case will be examined in a forthcoming paper \cite{MMV11}.

A few lines on the chosen context. We work internally to a semi-abelian category in which the
notions of Huq commutator and Smith commutator coincide.
This allows us, among other things, to use a simplified version of internal crossed modules
without loosing the equivalence with internal groupoids (see \cite{MFVdL10}).
The categories of  groups, Lie algebras, rings and many other algebraic structures not only in
Set, but
in any Grothendieck topos satisfy this condition (see Remark \ref{rem:HuqIsSmith}), so that
our context include also that of \cite{AN09}.

Finally, let us give a glance to possible developments of the present work.
Quite a lot of higher dimensional group theory has been developed starting from the pioneer
works of P. Deligne \cite{Del73} and A. Fr\"ohlich and C.T.C. Wall \cite{FW74}  on Picard
categories (also called 2-groups or
categorical groups), taking monoidal functors as morphisms (see for example
\cite{Vitale02}, \cite{AN09}, \cite{Elg11}
and the references therein). On the other hand, group theory has been the paradigmatic
example to
develop in recent years semi-abelian categorical algebra (see Section 9 and the
references therein).
The fact of disposing of an internal notion  of monoidal functor (the butterflies) should
make possible
to join these two generalizations of group theory and to develop a ``higher dimensional
semi-abelian
categorical algebra" which could cover as special cases most of the known results on
(strict) categorical
groups and (strict) Lie 2-algebras.

The layout of this paper is as follows: in Section 2 we recall the equivalence between
internal groupoids and internal crossed module, a result due to G. Janelidze
(see \cite{Janelidze03}) and which holds in any semi-abelian category; in Sections
3 and 4 we study
the bicategory $\cB(\cC)$ of butterflies in a semi-abelian category $\cC$ with
``Huq = Smith'';
in Section 5 we prove that $\cB(\cC)$ is the bicategory of fractions of internal functors
with respect to weak equivalences;
in section 6 we examine the three leading examples of groups,
rings and Lie algebras; section 7 is a short section devoted to the
classification of extensions which follows from Section 5; in Section 8 we specialize the main result of Section 5 to the case
where $\cC$ is a free exact category; finally, Section 9 is a reminder on protomodular and
semi-abelian
categories.
The reader who is not familiar with semi-abelian categories should have a glance to
Section 9 before reading Section 2.\\

\noindent {\em Notation:} the composite of
$f \colon A \to B$ and $g \colon B \to C$ is written $f \cdot g$ or $fg.$

\noindent {\em Terminology:} bicategory means bicategory with invertible 2-cells.

\section{Internal groupoids and internal crossed modules}

An introduction to internal categories can be found in Chapter 8 of \cite{Borceux94}.
For basic facts on 2-categories and bicategories, see \cite{Benabou67} or
Chapter 7 in \cite{Borceux94}.

\begin{IntGroupoid}\label{IntGroupoid}{\em
Let $\cC$ be a category with finite limits. We use the following notation:
\begin{enumerate}
\item[-] An (internal) groupoid $\G$ in $\cC$ is depicted as
$$\xymatrix{G_1 \times_{c,d}G_1 \ar[r]^>>>>>>m & G_1 \ar@<1ex>[r]^d
\ar@<-1ex>[r]_c & G_0 \ar[l]|e & G_1 \ar[r]^{i} & G_1}$$
where
$$\xymatrix{G_1 \times_{c,d}G_1 \ar[r]^>>>>>>{\pi_2} \ar[d]_{\pi_1}
& G_1 \ar[d]^d \\ G_1 \ar[r]_c & G_0}$$
is a pullback;
\item[-] An (internal) functor $P=(p_1,p_0) \colon \G \to \bH$ is depicted as
$$\xymatrix{G_1 \ar[r]^{p_1}
\ar@<-0.5ex>[d]_d \ar@<0.5ex>[d]^c & H_1 \ar@<-0.5ex>[d]_d \ar@<0.5ex>[d]^c \\
G_0 \ar[r]_{p_0} & H_0}$$
\item[-] An (internal) natural transformation
$\alpha \colon P \Rightarrow Q \colon \G \to \bH$ is depicted as
$$\xymatrix{G_1 \ar@<0.5ex>[r]^{p_1} \ar@<-0,5ex>[r]_{q_1}
\ar@<-0.5ex>[d]_d \ar@<0.5ex>[d]^c & H_1 \ar@<-0.5ex>[d]_d \ar@<0.5ex>[d]^c  \\
G_0 \ar@<0.5ex>[r]^{p_0} \ar@<-0,5ex>[r]_{q_0} \ar[ru]|{\alpha} & H_0}$$
\end{enumerate}
Internal groupoids,  functors and  natural transformations
form a 2-category (with invertible 2-cells) denoted by $\GrpdC.$
}\end{IntGroupoid}
When dealing with internal structures, it is sometimes useful to use {\em virtual}
objects and arrows as if those would be internal to the category of sets. For instance,
we could describe the object $G_1 \times_{c,d}G_1$  as the ``set''of composable arrows
 $\xymatrix{\cdot\ar[r]^f&\cdot\ar[r]^g&\cdot}$. Yoneda embedding
makes this precise, as explained in \cite{BB}, Metatheorem 0.2.7.

\begin{DiscFib}\label{DiscFib}{\em
An internal functor $P \colon \G \to \bH$ as above is called a {\em discrete cofibration}
when the commutative square $p_1 d = d p_0$ is a pullback. Dually, $P$ is a {\em discrete fibration}
when the square $p_1 c = c p_0$ is a pullback; for  groupoids these two notions
are equivalent.

Any  groupoid comes with several canonical fibrations onto it, some of those
being of interest in the rest of the paper.

For instance, let us consider the diagram
$$\xymatrix{
R[c] \ar[r]^{\widetilde{d}} \ar@<-0.5ex>[d]_{c_1} \ar@<0.5ex>[d]^{c_2}
& G_1 \ar@<-0.5ex>[d]_d \ar@<0.5ex>[d]^c \\
G_1\ar[r]_{d} & G_0
}$$
where $(R[c],c_1,c_2)$ is a kernel pair of $c$, and $\widetilde{d}$ is the morphism that sends
the pair of converging virtual arrows $\xymatrix{x\ar[r]^f&y&z\ar[l]_g}$ in the composition
$\xymatrix{x\ar[r]^{f\cdot g^{-1}}&z}$. The pair $(\widetilde{d},d)$ is clearly a discrete fibration
 of  groupoids. A similar argument can be developed for $R[d]$.

}\end{DiscFib}

\begin{WeakEquiv}\label{WealEquiv}{\em The following notion is due to M. Bunge
and R. Par\'e, see \cite{BP79}: An internal functor $P \colon \G \to \bH$ is a weak
equivalence if it is
\begin{enumerate}
\item[-] (internally) full and faithful, that is, the diagram
$$\xymatrix{ & & G_1 \ar[lld]_{d} \ar[d]^{p_1} \ar[rrd]^{c} \\
G_0 \ar[rd]_{p_0} & & H_1 \ar[ld]^{d} \ar[rd]_{c} & & G_0 \ar[ld]^{p_0} \\
 & H_0 & & H_0}$$
 is a limit, and
\item[-]  (internally) essentially surjective on objects, that is,
$$\xymatrix{G_0 \times_{p_0,c} H_1 \ar[r]^>>>>>{t_2} & H_1 \ar[r]^d & H_0}$$
is a regular epimorphism, where
$$\xymatrix{G_0 \times_{p_0,c} H_1 \ar[r]^>>>>>{t_2} \ar[d]_{t_1}
& H_1 \ar[d]^c \\
G_0 \ar[r]_{p_0} & H_0}$$
is a pullback.
\end{enumerate}
Observe that $P \colon \G \to \bH$ is (internally) full and faithful iff for every
 groupoid $\A$ the functor
$$- \cdot P \colon \GrpdC(\A,\G) \to \GrpdC(\A,\bH)$$
is full and faithful in the usual sense. Moreover, a weak equivalence $P$ is an
equivalence iff
$$\xymatrix{G_0 \times_{p_0,c} H_1 \ar[r]^>>>>>{t_2} & H_1 \ar[r]^d & H_0}$$
is a split epimorphism.
}\end{WeakEquiv}

From now on we assume that $\cC$ is a semi-abelian category in which the implication
``Huq $\Rightarrow$ Smith'' holds (for undefined notions and notations concerning
semi-abelian categories the reader is addressed to Section 9).

\begin{XMod}\label{XMod}{\em An (internal) crossed module $\G$ in $\cC$ is given
by a morphism $\partial \colon G \to G_0$ and an action
$\xi \colon G_0 \flat G \to G$ such that each part of the diagram
$$\xymatrix{G \flat G \ar[d]_{\partial \flat 1} \ar[rr]^{\chi_{{}_{G}}} & & G \ar[d]^1 \\
G_0 \flat G \ar[d]_{1 \flat \partial} \ar[rr]^{\xi} & & G \ar[d]^{\partial} \\
G_0 \flat G_0 \ar[rr]_{\chi_{{}_{G_0}}} & & G_0}$$
commutes, $\chi_{{}_{X}}$ being the canonical \emph{conjugation action} for the object $X$.
The commutativity of the upper part is called Peiffer condition, the
commutativity of the lower part is called \emph{precrossed module condition}.
\par
A morphism $P \colon \bH \to \G$ of  crossed modules is given by morphisms
$p \colon H \to G$ and $p_0 \colon H_0 \to G_0$ such that each part of the diagram
$$\xymatrix{H_0 \flat H \ar[d]_{\xi} \ar[rr]^{p_0 \flat p} & &
G_0 \flat G \ar[d]^{\xi} \\
H \ar[d]_{\partial} \ar[rr]^{p} & & G \ar[d]^{\partial} \\
H_0 \ar[rr]_{p_0} & & G_0}$$
commutes. In the following, we will refer to the upper commutative square above by saying that
 the pair $(p,p_0)$ is \emph{equivariant w.r.t.\ the actions}.
\par
Internal crossed modules with their morphisms form a category denoted by $\XModC.$
}\end{XMod}

\begin{Remark}\label{RemPeiffer}{\em In an arbitrary semi-abelian category the notion
of  crossed module introduced by G. Janelidze in \cite{Janelidze03} is stronger
than the one in \ref{XMod}. The notion we use (already considered in \cite{Janelidze03}
and further studied in \cite{MM10}) is equivalent to the original one thanks to the
condition ``Huq $\Rightarrow$ Smith'', as proved in \cite{MFVdL10}.
}\end{Remark}

In the following proposition we consider $\GrpdC$ as a category, that is, we forget
natural transformations.

\begin{Proposition}\label{PropGrpdXMod}(Janelidze \cite{Janelidze03})
The categories $\GrpdC$ and $\XModC$ are equivalent.
\end{Proposition}

\begin{proof} (\emph{sketch})  Let
$$\xymatrix{G_1 \times_{c,d}G_1 \ar[r]^>>>>>>m & G_1 \ar@<1ex>[r]^d
\ar@<-1ex>[r]_c & G_0 \ar[l]|e}$$
be a  groupoid and consider the following commutative diagram, where the
rows are kernel diagrams,
\begin{equation}\label{Diag:gpd_to_act}
\raisebox{4ex}{\xymatrix{G_0 \flat G \ar[rr]^{j_{G_0,G}} \ar[d]_{\exists !}^{\xi}
& & G_0 + G \ar[d]^{[e,g]} \ar[rr]^{[1,0]} & & G_0 \ar[d]^1 \\
G \ar[rr]_{g} & & G_1 \ar[rr]_c & & G_0}}\end{equation}
We obtain a  crossed module
$$\xymatrix{G \ar[r]^{g} & G_1 \ar[r]^{d} & G_0 &
G_0 \flat G \ar[r]^{\xi} & G}$$
This describes the equivalence functor
$$J \colon \GrpdC \to \XModC$$
on objects; its extension to arrows is straightforward.
\par
Conversely, let
$$\xymatrix{G_0 \flat G \ar[r]^{\xi} & G \ar[r]^{\partial} & G_0}$$
be a  crossed module and consider the semi-direct product, given by the coequalizer
$$\xymatrix{G_0 \flat G \ar[rr]^{j_{G_0,G}} \ar[rd]_{\xi} & & G_0 + G
\ar[r]^{q_{\xi}} & G \rtimes_{\xi} G_0 \\
& G \ar[ru]_{i_G} }$$
We obtain a  reflexive graph
$$G_1=\xymatrix{G \rtimes_{\xi} G_0 \ar@<1ex>[rr]^d \ar@<-1ex>[rr]_c & &
G_0 \ar[ll]|e}$$
where $c$ and $e$ are, respectively, the canonical projection from and
the canonical injection into the semi-direct product, and $d$
is the unique morphism such that the diagram
$$\xymatrix{G \ar[rr]^{g} \ar[rrd]_{\partial} & & G \rtimes_{\xi} G_0
\ar[d]_{d} & & G_0 \ar[ll]_{e} \ar[lld]^1 \\
& & G_0}$$
commutes. For a detailed proof, see \cite{Janelidze03}.
\end{proof}
We will refer to the functor $J$ as to the {\em normalization functor},
and to its quasi-inverse as to the {\em denormalization functor}.
\begin{Notation}{\em
Here and in the following we will denote kernel of the codomain arrows with the lower case letter of the groupoid involved,
e.g.\ the following sequence is exact:
$$
\xymatrix{G\ar[r]^g&G_1\ar[r]^c&G_0}
$$
 Moreover $g\cdot i=\ker d$  will be often denoted $g^{\bullet}$.}
\end{Notation}
\begin{2CatXMods}\label{2CatXMods}{\em
The category $\XModC$ has an obvious 2-categorical structure. In fact it suffices to translate the notion
of natural transformation for internal functors in the language of crossed module in order to obtain
the 2-cells of $\XModC$.
\begin{Definition}\label{Def2cellXmod}
Consider two parallel morphisms $P, Q \colon \bH \rightrightarrows \G$
of  crossed modules.
An arrow $\alpha \colon H_0 \to G_1=G\rtimes_{\xi} G_0$ is a  natural transformation
between $P$ and $Q$ if
$\alpha \cdot d = p_0,$ $\alpha \cdot c = q_0,$ and the diagram
$$\xymatrix{H \ar[r]^{p}\ar@{}[dr]|{(\star)}
\ar[d]_{\langle \partial \cdot \alpha, q \cdot g \rangle} & G\ar[d]^{g} \\
G_1 \times_{c,d} G_1 \ar[r]_>>>>>>{m} & G_1}$$
commutes.
\end{Definition}
The definition above involves several
non-elementary constructions, such as the semidirect product $G_1=G\rtimes G_0$ and the morphisms $d$ and $m$.
That is why we give an equivalent but easier to handle version of the diagram $(\star)$. This is done by the following
\begin{Proposition}\label{PropNatPFFGraph}
Diagram $(\star)$ of Definition \ref{Def2cellXmod} commutes iff the following one does:
$$
\xymatrix{
H\ar[r]^{\partial}\ar[d]_{\langle p,q\rangle}
&H_0\ar[d]^{\alpha}
\\
G\times G\ar[r]_{m_0}&G_1
}
$$
where $m_0= g \sharp g^\bullet$ is the cooperator of the maps $\xymatrix{G\ar[r]^{g}&G_1&\ar[l]_(.4){g^{\bullet}}G}$
{\em (see section \ref{HuqIsSmith})}.
\end{Proposition}
For the proof, it suffices to compute with elements and then apply Yoneda embedding.
To this end, recall that $m_0$ is the morphism that sends the pair of arrows
$(\xymatrix{x\ar[r]^a&0},\xymatrix{y\ar[r]^b&0})$ to the composition $\xymatrix@C=6ex{x\ar[r]^{a \cdot b^{-1}}&y}$, in other words
$m_0$ performs a division $a/b$. Together with the conditions on domain and codomain,
 the last commutative diagram defines a  transformation of {\em Peiffer graphs} \cite{MM10}.

\begin{Lemma}\label{LemmaTrasfNat}
Arrows satisfying Definition \ref{Def2cellXmod} correspond biunivocally to  natural transformations between
the internal functors determined by the morphisms $P$ and $Q$.
\end{Lemma}
\begin{proof}
Recall (for instance, from \cite{Borceux94}) that a  natural transformation between two internal
functors $P=(p_1,p_0)$ and $Q=(q_1,q_0)\colon$ $\bH\to\G $ is defined as a
morphism $\alpha\colon H_0 \to G_1$ satisfying $\alpha \cdot d = p_0,$ $ \alpha \cdot c = q_0,$ and
such that the diagram
$$
\xymatrix@C=10ex{\ar@{}[dr]|{(\star\star)}
H_1 \ar[r]^{\langle p_1, c\cdot \alpha \rangle}  \ar[d]_{\langle d\cdot\alpha,q_1\rangle}
&G_1 \times_{c,d} G_1 \ar[d]^m
\\
G_1 \times_{c,d} G_1\ar[r]_m
&G_1
}
$$
commutes. So we have to prove that $(\star)$ commutes iff $(\star\star)$ commutes. The ``if'' part is dealt with
by simply precomposing the diagram $(\star\star)$ with the monomorphism $h\colon H\to H_1$, conversely one observes
that since the base category is protomodular, the pair $$h,e\colon H_0\to H_1$$ is (strongly) jointly epic, so that
$(\star\star)$ commutes iff it commutes when precomposed both with $h$ and $e$. The first precomposition
is precisely $(\star)$, the second one is trivial.
\end{proof}
In conclusion we have proved the following
\begin{Corollary}\label{Grd_bieq_Xmod}
The equivalence between $\GrpdC$ and $\XModC$ extends to a biequivalence.
\end{Corollary}
}\end{2CatXMods}

\begin{Remark}\label{RemWeakEquiv}\label{RemDiscrFibr}{\em
\begin{itemize}
\item[(i)] It has been proved in \cite{EKVdL05} that,
in the equivalence $\GrpdC \simeq \XModC,$ a morphism
$P \colon \bH \to \G$ of  crossed modules corresponds to a weak equivalence
iff the arrows induced on kernels and cokernels are isomorphisms:
$$\xymatrix{\ker \partial \ar[r]^{\simeq} \ar[d] & \ker \partial \ar[d] \\
H \ar[r]^{p} \ar[d]_{\partial} & G \ar[d]^{\partial} \\
H_0 \ar[r]^{p_0} \ar[d] & G_0 \ar[d] \\
\coker \partial \ar[r]^{\simeq} & \coker \partial}$$

\item[(ii)] Under the mentioned equivalence, it is easy to show that a morphism
$P \colon \bH \to \G$ of  crossed modules corresponds to a discrete fibration
iff $p \colon H \to G$ is an isomorphism.
\end{itemize}
}\end{Remark}

\begin{Lemma}\label{LemmaPBXmod}
Let $(\partial \colon H \to H_0, \;\xi \colon H_0 \flat H \to H)$ be a  crossed module and consider a morphism $\sigma \colon E \to H.$
Consider also the pullback
$$\xymatrix{E \times_{\sigma,\partial}H \ar[rr]^{\overline\sigma}
\ar[d]_{\overline\partial} & & H \ar[d]^{\partial} \\
E \ar[rr]_{\sigma} \ar@{}[rru]|{(*)} & & H_0}$$
\begin{enumerate}
\item The morphism
$\overline\partial \colon E \times_{\sigma,\partial}H \to E$
can be equipped with a canonical action
$$\overline\xi \colon E \flat (E \times_{\sigma,\partial}H)
\to E \times_{\sigma,\partial}H$$
in such a way that the pair $(\overline\partial, \,\overline\xi\,)$
is a  crossed module and the diagram $(*)$ is a morphism of
crossed modules.
\item Moreover, if $\sigma \colon E \to H_0$ is a regular epimorphism, then $(*)$
is a weak equivalence.
\end{enumerate}
\end{Lemma}

\begin{proof}
1. By naturality of $\chi$ and the precrossed module condition on
$(\partial, \xi),$ the diagram
$$\xymatrix{E \flat (E \times_{\sigma,\partial}H)
\ar[r]^>>>>>{\sigma \flat \overline\sigma} \ar[d]_{1 \flat \overline\partial}
& H_0 \flat H \ar[r]^{\xi} & H \ar[d]^{\partial} \\
E \flat E \ar[r]_{\chi_E} & E \ar[r]_{\sigma} & H_0}$$
commutes. Therefore, by the universal property of the pullback, we get a unique morphism
$$\overline\xi \colon E \flat (E \times_{\sigma,\partial}H)
\to E \times_{\sigma,\partial}H$$
making  the following diagram commute
$$\xymatrix{E \flat E \ar[d]_{\chi_E} &
E \flat (E \times_{\sigma,\partial}H) \ar[d]_{\overline\xi}
\ar[l]_>>>>>{1 \flat \overline\partial}
\ar[r]^>>>>>{\sigma \flat \overline\sigma}
& H_0 \flat H \ar[d]^{\xi} \\
E \ar@{}[ru]|{(1)} & E \times_{\sigma,\partial}H
\ar[l]^>>>>>>>>>{\overline\partial} \ar@{}[ru]|{(2)}
\ar[r]_>>>>>>>>>{\overline\sigma} & H\ .}$$
Since commutativity of $(1)$ is the precrossed module condition on
$(\overline\partial, \,\overline\xi\,)$ and
commutativity of $(2)$ says that $(*)$
is a morphism of crossed modules, it remains to check Peiffer condition on
$(\overline\partial, \,\overline\xi\,),$ i.e., the commutativity of
$$\xymatrix{(E \times_{\sigma,\partial}H) \flat (E \times_{\sigma,\partial}H)
\ar[rr]^>>>>>>>>>>{\chi}
\ar[d]_{\overline\partial \flat 1}
& & E \times_{\sigma,\partial}H \ar[d]^{1} \\
E \flat (E \times_{\sigma,\partial}H) \ar[rr]_{\overline\xi}
& & E \times_{\sigma,\partial}H\ .}$$
For this, it suffices to compose with the pullback projections
$$\xymatrix{E & E \times_{\sigma,\partial}H \ar[l]_>>>>>>{\overline\partial}
\ar[r]^>>>>>>{\overline\sigma} & H\ .}$$
The diagram obtained by composing with $\overline\partial$ commutes by naturality
of $\chi$ and commutativity of $(1),$ the one obtained by composing with
$\overline\sigma$ commutes by naturality of $\chi,$ commutativity of $(2)$
and Peiffer condition on $(\partial, \xi).$ \\
2. Since $(*)$ is a pullback, kernels of parallel arrows are isomorphic. Henceforth, since
$\sigma$ and $\overline\sigma$ are regular epimorphisms, $(*)$ is also a pushout, so that
the induced arrow $\coker \overline\partial \to \coker \partial$ is an isomorphism.
From \ref{RemWeakEquiv} we conclude that $(*)$ is a weak equivalence.
\end{proof}

\section{The bicategory of butterflies}

We are going to describe the bicategory
$\cB(\cC)$ of  crossed modules and  butterflies in $\cC.$

\begin{Farfalle}\label{Farfalle}{\em The following notion has been introduced,
when $\cC$ is the category of groups, by B. Noohi in \cite{Noohi05}, see also
\cite{AN09} (a special case of butterflies was used by D. F. Holt in
\cite{Holt79} to classify group extensions.). Let $\G$ and $\bH$ be  crossed modules.
A (internal) butterfly from $\bH$ to $\G$ is given by a commutative diagram of the form
$$\xymatrix@!=2ex{H \ar[rd]^{\kappa} \ar[dd]_{\partial} & &
G \ar[dd]^{\partial} \ar[ld]_{\iota} \\
& E \ar[ld]^{\sigma} \ar[rd]_{\rho} \\
H_0 & & G_0}$$
such that
\begin{enumerate}
\item[i.] $\kappa \cdot \rho = 0$, \ \ \ \ i.e.\ $(\kappa, \rho)$ is a complex
\item[ii.] $\iota = \ker \sigma$ and $\sigma = \coker \iota$,  \ \ \ \ i.e.\ $(\iota, \sigma)$ is an extension
\item[iii.] the diagram
$$
\xymatrix{
E \flat H \ar[r]^{\sigma \flat 1}\ar[d]_{1\flat \kappa} & H_0 \flat H \ar[r]^{\xi} & H\ar[d]^{\kappa}\\
E\flat E\ar[rr]_{\chi_E}&&E
}
$$
commutes, \ \ \ \ i.e.\ the pair $(\kappa,\sigma \flat 1\cdot\xi)$ is a precrossed module,
\item[iv.]  the diagram
$$
\xymatrix{
E \flat G \ar[r]^{\rho \flat 1}\ar[d]_{1\flat \iota} & G_0 \flat G \ar[r]^{\xi} & G\ar[d]^{\iota}\\
E\flat E\ar[rr]_{\chi_E}&&E}
$$
commutes, \ \ \ \ i.e. \ the pair $(\iota,\rho \flat 1\cdot \xi)$ is a precrossed module.
\end{enumerate}
When no confusion arises, we denote a butterfly
$(E,\kappa, \rho, \iota, \sigma)$ from $\bH$ to $\G$ simply by
$$E \colon \bH \to \G$$
Observe that since $(\iota,\sigma)$ is an exact sequence, the situation is not as symmetrical as
it may appear at first sight. Actually  $\iota$ is a mono, then the action $\rho \flat 1\cdot \xi$
is nothing but the conjugation action $\chi_E$ restricted to the subobject $G$. Moreover
$\iota$ can be recovered as the normalization of the
equivalence relation $(R[\sigma],\sigma_1,\sigma_2)$, i.e. $\iota = \ker(\sigma_2)\cdot \sigma_1$.

Given butterflies $E, E' \colon \bH \to \G,$ a morphism of butterflies is a morphism
$$f \colon E \to E'$$
such that the diagrams
$$\xymatrix{& E\ar[dd]^f \ar[rd]^{\rho} \\
H \ar[ru]^{\kappa} \ar[rd]_{\kappa'} & & G_0 \\
& E' \ar[ru]_{\rho'} }
\;\;\;\;\;
\xymatrix{& E\ar[dd]^f \ar[ld]_{\sigma} \\
H_0 & & G \ar[lu]_{\iota} \ar[ld]^{\iota'} \\
& E' \ar[lu]^{\sigma'} }$$
commute. In particular, $f$ is a morphism of extensions, so that by the short five lemma
it is an isomorphism.
}\end{Farfalle}

\begin{Remark}\label{RemEquivarianzaFarf}{\em
Conditions (iii) and (iv) in \ref{Farfalle}  imply that the pairs
$(\kappa, \sigma \flat 1 \cdot \xi)$ and $(\iota, \rho \flat 1\cdot \xi)$
are indeed  crossed modules and, therefore,
$$
\xymatrix{H \ar[d]_{\partial} & H \ar[d]^{\kappa} \ar[l]_1 \\
H_0 & E \ar[l]^{\sigma}}
\; \raisebox{-3.5ex}{ and } \;\;
\xymatrix{G \ar[d]_{\iota} \ar[r]^1 & G \ar[d]^{\partial} \\
E \ar[r]_{\rho} & G_0}
$$
are morphisms of  crossed modules, hence by \ref{RemWeakEquiv} discrete fibrations.
}\end{Remark}

\begin{Torsori}\label{Torsori}{\em
Using the equivalence between  crossed modules and  groupoids
described in \ref{PropGrpdXMod}, butterflies correspond to {\em fractors}, i.e.\
diagrams of the form
$$\xymatrix@!{ & R \ar[ld]_{\overline \sigma} \ar@<-0,5ex>[rd]_{c}
\ar@<0,5ex>[rd]^{d} & & R[\sigma] \ar@<-0,5ex>[ld]_{\sigma_1}
\ar@<0,5ex>[ld]^{\sigma_2} \ar[rd]^{\overline \rho} \\
H_1 \ar@<-0,5ex>[rd]_{c} \ar@<0,5ex>[rd]^{d} & & E \ar[ld]^{\sigma}
\ar[rd]_{\rho} & & G_1 \ar@<-0,5ex>[ld]_{d} \ar@<0,5ex>[ld]^{c} \\
& H_0 & & G_0}$$
where
\begin{enumerate}
\item  functors $(\overline{\sigma},\sigma)$ and $(\overline{\rho},\rho)$
are discrete fibrations,
\item $\sigma$ is a regular epimorphism, and $R[\sigma]$ is its kernel pair,
\item $\rho$ coequalizes
$d, c \colon R \rightrightarrows E.$
\end{enumerate}
}\end{Torsori}

\begin{proof}
Denormalizing the morphisms of crossed modules in the diagram
$$
\xymatrix{
&H \ar[dl]_1\ar[dr]^{\kappa}
&&G\ar[dl]_{\iota}\ar[dr]^1
\\
H\ar[dr]_{\partial}
&&E\ar[dl]^{\sigma}\ar[dr]_{\rho}
&&G\ar[dl]^{\partial}
\\
&H_0&&G_0
}
$$
one easily gets a fractor as above, where
$H_1=H \rtimes_{\xi}H_0$, $G_1=G \rtimes_{\xi}G_0$ and $R=H \rtimes_{\sigma \flat 1\cdot \xi}E$.
The fact that the groupoid associated to $\iota$ is isomorphic to $(R[\sigma],\sigma_1,\sigma_2)$
is due to the fact that $\iota=\ker\sigma$.
Finally, $\rho$ coequalizes $d$ and $c$ since the pair
$$\xymatrix{H\ar[r]^{\langle h,0\rangle}&R&E\ar[l]_{e}}$$ is jointly (strongly) epic,
by protomodularity.

Conversely, starting from a fractor as above, we get the butterfly
$$\xymatrix{H \ar[dd]_{\partial} \ar[rrd]^{\langle h,0 \rangle \cdot d} & & &
\ker \sigma \simeq G \ar[ld]_{j} \ar[dd]^{\partial} \\
& & E \ar[lld]^{\sigma} \ar[rd]_{\rho} \\
H_0 & & & G_0}$$
where $\langle h,0 \rangle \colon H \to R$ comes from the universal property
of the pullback
$$\xymatrix{H_1 \ar[d]_{c} & R \ar[l]_{\overline \sigma} \ar[d]^{c} \\
H_0 & E \ar[l]^{\sigma}}$$
and the isomorphism $\ker \sigma \simeq G$ is the composite of the
following isomorphisms determined by bottom pullback squares:
$$\xymatrix{\ker \sigma \ar[d]_{j} & \ker \sigma_2 \ar[d] \ar[l]_{\simeq}
\ar[r]^{\simeq} & G \ar[d]^{g} \\
E\ar[d]_{\sigma} & R[\sigma]\ar[d]^{\sigma_2} \ar[l]_{\sigma_1} \ar[r]^{\overline \rho} & G_1\ar[d]^{c}
\\
H_0 & E\ar[l]^{\sigma}\ar[r]_{\rho}&G_0
}$$
\end{proof}

\begin{Remark}
In a recent paper by D. Bourn \cite{Bourn10}, what we have called fractor is termed {\em left regularly faithful profunctor}
\end{Remark}

\begin{Remark}\label{RemProfunctors}{\em
Given a fractor as in \ref{Torsori}, one can consider also the kernel pair of the map
$\overline \sigma$ and perform the construction below, where dashed arrows are suitably obtained by the
universal property of the kernel pair $R[\sigma]$,
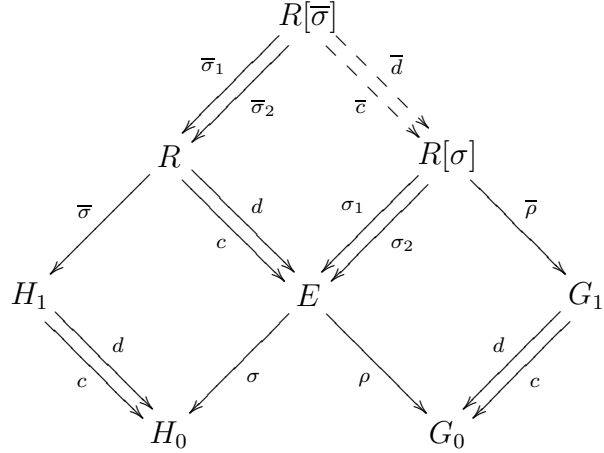
$$\xymatrix@!{
&&R[\overline \sigma]\ar@<-0,5ex>[ld]_{\overline{\sigma}_1}\ar@<0,5ex>[ld]^{\overline{\sigma}_2}
\ar@{-->}@<-0,5ex>[rd]_{\overline{c}}\ar@{-->}@<0,5ex>[rd]^{\overline{d}}\\
 & R \ar[ld]_{\overline \sigma} \ar@<-0,5ex>[rd]_{c}
\ar@<0,5ex>[rd]^{d} & & R[\sigma] \ar@<-0,5ex>[ld]_{\sigma_1}\ar@<0,5ex>[ld]^{\sigma_2} \ar[rd]^{\overline \rho} \\
H_1 \ar@<-0,5ex>[rd]_{c} \ar@<0,5ex>[rd]^{d} & & E \ar[ld]^{\sigma}
\ar[rd]_{\rho} & & G_1 \ar@<-0,5ex>[ld]_{d} \ar@<0,5ex>[ld]^{c} \\
& H_0 & & G_0}$$
One finds out that the central square is a double groupoid over $E$. More precisely, it is a
{\em centralizing double groupoid}, as defined by D. Bourn in \cite{Bourn10}, since
$(\overline{\sigma},\sigma)$ is a discrete fibration. Together with the two other squares,
this gives rise to a particular profunctor $\bH \looparrowright \G$ of  groupoids
(profunctors were introduced by J. B\'enabou with the name of {\it distributeurs} \cite{Benabou73},
the internal version  can be found in \cite{PtJ}).
The precise relationship between butterflies and profunctors will be described in a forthcoming paper
\cite{MMV11}.
}
\end{Remark}
\begin{IdFarfalle}\label{IdFarfalle}{\em
We are going to prove that the canonical fractor associated to a  groupoid gives  the
identity butterfly associated to a  crossed module.

Let $\G$ be a  crossed module. In order to construct the identity butterfly on $\G,$ consider the  groupoid
associated to $\G$ as in the proof of Proposition \ref{PropGrpdXMod}
$$G_1=\xymatrix{G \rtimes_{\xi} G_0 \ar@<1ex>[rr]^d \ar@<-1ex>[rr]_c & &
G_0 \ar[ll]|e}$$
Then, following the example described in \ref{DiscFib}, one can associate  $\G$ to the fractor
$$
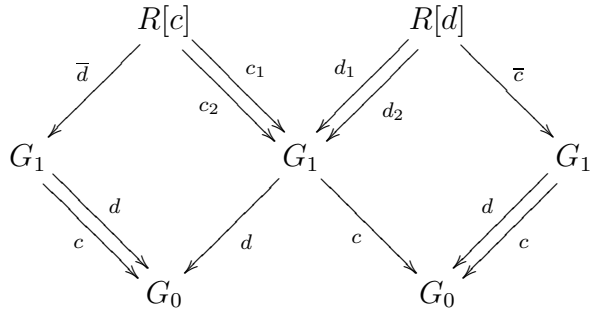
\xymatrix@!{ & R[c] \ar[ld]_{\overline{d}} \ar@<-0,5ex>[rd]_{c_2}
\ar@<0,5ex>[rd]^{c_1} & & R[d] \ar@<-0,5ex>[ld]_{d_1}
\ar@<0,5ex>[ld]^{d_2} \ar[rd]^{\overline{c}} \\
G_1 \ar@<-0,5ex>[rd]_{c} \ar@<0,5ex>[rd]^{d} & & G_1 \ar[ld]^{d}
\ar[rd]_{c} & & G_1 \ar@<-0,5ex>[ld]_{d} \ar@<0,5ex>[ld]^{c} \\
& G_0 & & G_0}$$
The butterfly associated to this fractor, by means of the normalization process described in
\ref{Torsori}, is called the identity butterfly of the crossed module $\G$.
Actually it acts as an identity w.r.t.\ the composition described in \ref{CompFarfalle}.
It is represented explicitly in the diagram below:
$$\xymatrix{G \ar[dd]_{\partial} \ar[rd]^{g} & &
G \ar[dd]^{\partial} \ar[ld]_{g^\bullet} \\
& G_1 \ar[ld]^d \ar[rd]_c \\
G_0 & & G_0}$$

Actually, in this paper, we will use as identity butterfly the isomorphic
contravariant version of the one above:
$$\xymatrix{G \ar[dd]_{\partial} \ar[rd]^{g^\bullet} & &
G \ar[dd]^{\partial} \ar[ld]_{g} \\
& G_1 \ar[ld]^c \ar[rd]_d \\
G_0 & & G_0}$$
the isomorphism being realized by the inverse map $i\colon G_1\to G_1$.
This choice does not affect the computations (consider that the composition will be defined only up to
isomorphisms), but is coherent with the  normalization of a groupoid via the kernel of the codomain.

}\end{IdFarfalle}

\begin{CompFarfalle}\label{CompFarfalle}{\em Let
$$\xymatrix{\bH \ar[r]^{E} & \G \ar[r]^{E'} & \bK}$$
be butterflies. In order to construct their composition, consider the diagram
$$\xymatrix{& & Q
\ar@/^6pc/[rrdddd]^{\overline{s \rho'}} \ar@/_6pc/[lldddd]_{\overline{r \sigma}} \\
& & E \times_{\rho, \sigma'} E' \ar[u]_q \ar[ldd]_r \ar[rdd]^s \\
H \ar[rru]^{\langle \kappa,0 \rangle} \ar[dd]_{\partial} \ar[rd]^{\kappa} & &
G \ar[u]|>>>>>>{\langle \iota,\kappa' \rangle}
\ar[ld]^{\iota} \ar[rd]_{\kappa'} \ar[dd]_{\partial} & &
K \ar[llu]_{\langle 0,\iota' \rangle} \ar[ld]_{\iota'} \ar[dd]^{\partial} \\
& E \ar[ld]^{\sigma} \ar[rd]_{\rho} & & E' \ar[ld]^{\sigma'} \ar[rd]_{\rho'} \\
H_0 & & G_0 & & K_0}$$
where
\begin{enumerate}
\item[-] $E \times_{\rho, \sigma'} E'$ is the pullback of $\rho$ and $\sigma',$
with projections $r$ and $s$, so that
$$
\langle \kappa,0 \rangle \cdot s = 0 = \langle 0,\iota' \rangle \cdot r, \;
\langle \kappa,0 \rangle \cdot r = \kappa, \;
\langle 0,\iota' \rangle \cdot s =\iota', \;
\langle \iota,\kappa' \rangle \cdot r = \iota, \;
\langle\iota,\kappa' \rangle \cdot s = \kappa',
$$
\item[-] $(Q,q)$ is the cokernel of $\langle\iota,\kappa'\rangle$,
\item[-] $\overline{r \sigma}$ and $\overline{s \rho'}$ are defined by
$q \cdot \overline{r \sigma} = r \cdot \sigma, \;
q \cdot \overline{s \rho'} = s \cdot \rho'.$
\end{enumerate}
The composition of $E$ and $E'$ is the butterfly
$$\xymatrix{H \ar[dd]_{\partial} \ar[rrd]^{\langle \kappa,0 \rangle \cdot q}
& & & & K \ar[lld]_{\langle 0,\iota' \rangle \cdot q} \ar[dd]^{\partial} \\
& & Q \ar[lld]^{\overline{r \sigma}}
\ar[rrd]_{\overline{s \rho'}} \\
H_0 & & & & K_0}$$
}\end{CompFarfalle}

\begin{proof} We have to check that the previous diagram is indeed a butterfly
from $\bH$ to $\bK.$
Commutativity of wings and condition \ref{Farfalle}.i are easy to check. \\
Condition \ref{Farfalle}.ii: first observe that
$\overline{r \sigma}$ is a regular epimorphism (because $\sigma$ and $\sigma'$
are), so that it is enough to show that
$\langle 0,\iota' \rangle \cdot q$ is the kernel of $\overline{r \sigma}.$
Since $\iota'$ is the kernel of $\sigma'$ and $r$ is a pullback of $\sigma',$ clearly
$\langle 0, \iota' \rangle$ is the kernel of $r.$
Consider now the commutative diagram
$$\xymatrix{G \ar[d]_{\iota} & G \ar[l]_{1}
\ar[d]^{\langle \iota, \kappa' \rangle} \\
E \ar[d]_{\sigma} \ar@{}[rd]|{(\star)}
& E \times_{\rho,\sigma'} E' \ar[d]^{q} \ar[l]_<<<<<{r} &
K = \ker r \ar[l]_<<<<<{\langle 0,\iota' \rangle} \ar[d]^{f} \\
H_0 & Q \ar[l]^{\overline{r \sigma}} &
\ker(\overline{r \sigma}) \ar[l]}$$
where $f$ is induced by the universal property of $\ker (\overline{r \sigma}).$
It suffices  to prove that $(\star)$ is a pullback. Indeed, if
$(\star)$ is a pullback, then $f$ is an isomorphism and, therefore,
$\langle 0,\iota' \rangle \cdot q$ is the kernel of $\overline{r \sigma}.$
Since $q$ and $\sigma$ are regular epimorphisms and $\iota$ is the kernel of $\sigma,$
 to show that $(\star)$ is a pullback
is equivalent to show that $\langle \iota,\kappa' \rangle$ is the kernel of $q.$
Since $\langle \iota,\kappa' \rangle$ is a monomorphism (because $\iota$ is),
 to prove that
$\langle \iota,\kappa' \rangle$ is a kernel (of its cokernel $q$) is equivalent to proving
that $\langle \iota,\kappa' \rangle$ is closed under conjugation in $E \times_{\rho,\sigma'}E'$
(see \cite{JMU07,MM10}).
The action of $E \times_{\rho, \sigma'} E'$ on $G$ is given by
$$\xymatrix{(E \times_{\rho,\sigma'} E') \flat G \ar[r]^>>>>>{r \flat 1} &
E \flat G \ar[r]^>>>>>{\rho \flat 1} & G_0 \flat G \ar[r]^>>>>>{\xi} & G}$$
or, equivalently, by
$$\xymatrix{(E \times_{\rho,\sigma'} E') \flat G \ar[r]^>>>>>{s \flat 1} &
E' \flat G \ar[r]^>>>>>{\sigma' \flat 1} & G_0 \flat G \ar[r]^>>>>>{\xi} & G}$$
Therefore, the normality of  $\langle \iota,\kappa' \rangle$ in
$E \times_{\rho,\sigma'}E'$ amounts to the commutativity of
$$\xymatrix{(E \times_{\rho,\sigma'} E') \flat G \ar[r]^>>>>>>>>>>{r \flat 1}
\ar[d]_{1 \flat \langle \iota, \kappa' \rangle} &
E \flat G \ar[r]^{\rho \flat 1} & G_0 \flat G \ar[r]^{\xi} &
G \ar[d]^{\langle \iota, \kappa' \rangle} \\
(E \times_{\rho,\sigma'} E') \flat (E \times_{\rho,\sigma'} E')
\ar[rrr]_>>>>>>>>>>>>>>>>>>>>>>>>{\chi}
& & & E \times_{\rho,\sigma'} E'}$$
For this, compose with the pullback projections $r$ and $s$ and use the naturality of
$\chi$ and, respectively, condition \ref{Farfalle}.iii on $\iota$ and
condition \ref{Farfalle}.iv on $\kappa'.$ \\
Condition \ref{Farfalle}.iii: since
$$q \flat 1 \colon (E \times_{\rho,\sigma'} E')
\flat H \to Q \flat H$$
is a (regular) epimorphism (see \cite{MM10}), condition \ref{Farfalle}.iii follows from
the commutativity of the whole diagram below
$$\xymatrix{ (E \times_{\rho,\sigma'} E') \flat H
\ar[r]^{q \flat 1} \ar[d]_{1 \flat \langle \kappa, 0 \rangle} &
Q  \flat H
\ar[r]^>>>>>{\overline{r \sigma} \flat 1} &
H_0 \flat H \ar[r]^{\xi} & H \ar[d]^{\langle \kappa, 0 \rangle} \\
(E \times_{\rho,\sigma'} E') \flat (E \times_{\rho,\sigma'} E')
\ar[rrr]_{\chi} \ar[d]_{q \flat q}
& & & E \times_{\rho,\sigma'} E' \ar[d]^{q} \\
Q \flat Q
\ar[rrr]_{\chi}
& & & Q}$$
The lower region commutes by naturality of $\chi.$ For the commutativity of the upper
region, compose with the pullback projections: composed with $s,$ both paths go
to zero; as far as $r$ is concerned, use condition \ref{Farfalle}.iii on $\kappa.$ \\
Condition \ref{Farfalle}.iv: same argument as for \ref{Farfalle}.iii.
\end{proof}

\begin{Proposition}\label{PropBicatFarfalle} We have a bicategory
$$\cB(\cC)$$
with internal crossed modules as objects,  butterflies as 1-cells, and morphisms
of butterflies as 2-cells.
\end{Proposition}

\begin{proof} Composition of butterflies and identity butterflies have been described
in \ref{CompFarfalle} and \ref{IdFarfalle}.
The rest of the proof is long but straightforward.
\end{proof}

Observe that in the identity butterfly (\ref{IdFarfalle}) both diagonals are extensions.
Butterflies with this property are called {\em flippable} (see \cite{Noohi05}).
\begin{Proposition}\label{PropEquivFarfalle} A flippable butterfly $E \colon \bH \to \G$
is an equivalence in the bicategory $\cB(\cC)$. A quasi-inverse $E^* \colon \G \to \bH$ is obtained
by twisting the wings of $E.$
\end{Proposition}

\begin{proof} Keep in mind \ref{CompFarfalle} and \ref{IdFarfalle}
and consider the diagram
$$\xymatrix{ & R[\rho] \ar[ldd]_{r} \ar[rdd]^{s} \\
& G \ar[u]|>>>>>>{\langle \iota, \iota \rangle} \ar[ld]^{\iota} \ar[rd]_{\iota} \\
E \ar[rd]_{\rho} & & E \ar[ld]^{\rho} \\
& G_0}$$
We have to prove that $H_1 = H \rtimes_{\xi} H_0$ is a cokernel of
$\langle \iota, \iota \rangle.$ Since $(\kappa, \rho)$ is an extension,
we can take (isomorphic) kernels in the left discrete fibration in the fractor corresponding to $E$
(see \ref{Torsori})
$$
\xymatrix{
H_1\ar@<-0,5ex>[d]_{c}\ar@<0,5ex>[d]^{d}
&R[\rho]\ar[l]_{\overline{\sigma}}\ar@<-0,5ex>[d]_{\rho_2}\ar@<0,5ex>[d]^{\rho_1}
&G\ar[l]_{\langle \iota,\iota \rangle}\ar[d]^1\\
H_0&E\ar[l]^{\sigma}&G\ar[l]^{\iota}
}
$$
Moreover, since $\sigma$ and its pullback $\overline{\sigma}$ are regular epimorphisms,
the horizontal rows are exact, and this concludes the proof.
\end{proof}

\section{Butterflies and morphisms of crossed modules}

In order to prove that $\cB(\cC)$ is the bicategory of fractions of $\GrpdC$
with respect to weak equivalences (Theorem \ref{ThMain}), we have to construct a
homomorphism of bicategories
$$\cF \colon \GrpdC \to \cB(\cC)$$
This task will be completed only in section \ref{HomF}, since before we have
 to provide some necessary constructions.

A preliminary step consists in associating a split butterfly (Definition \ref{DefSplitFarf})
to any morphism of  crossed modules.

\begin{Definition}\label{DefSplitFarf}{\em A butterfly $E \colon \bH \to \G$
is {\em split} when the extension
$$\xymatrix{H_0 & E \ar[l]_(.4){\sigma} & G \ar[l]_(.4){\iota}}$$
is split, that is, when there exists $s \colon H_0 \to E$ such that
$s \cdot \sigma = 1_{H_0}.$

A morphism of split butterfly is  simply a morphism of butterflies, so that it
need not  commute with sections.
}\end{Definition}

\begin{SplitFarfalle}\label{SplitFarfalle}{\em
Let $P \colon \bH \to \G$ be a morphism of  crossed modules.
We are going to construct a split butterfly $E_P \colon \bH \to \G.$

Consider the pullback:
$$\xymatrix{E_P \ar[r]^{\overline{p}} \ar[d]_{\sigma_P} & G_1 \ar[d]^{c} \\
H_0 \ar[r]_{p_0} & G_0}$$
If $\xi\colon G_0\flat G\to G$ is the action corresponding to the split epi $c\colon G_1\to G_0$,
it is easy to show that
$$\xymatrix{H_0 \flat G \ar[r]^{p_0 \flat 1} & G_0 \flat G \ar[r]^{\xi} & G}$$
is the action corresponding to the split epi $\sigma_P\colon E_P\to H_0$.
We get the split butterfly $E_P \colon \bH \to \G \colon$
$$
\xymatrix@C=16ex{\ar[dd]_{\partial}
H\ar[dr]^{\langle \partial ,p g i \rangle}
&&G\ar[dl]_{\langle 0,g \rangle}
\ar[dd]^{\partial}
\\
&E_P\ar[dl]^{\sigma_P} \ar[dr]_{\overline{p}d}
\\
H_0&&G_0}
$$
}\end{SplitFarfalle}

\begin{proof} Commutativity of the wings is  given by composing with pullback projections.\\
Condition \ref{Farfalle}.i: similarly one computes
$\langle \partial ,p g i \rangle\, \overline{p}\,d =   p\, g\, i\, d = p\, g\, c = p\, 0=0.$\\
Condition \ref{Farfalle}.ii: the North East - South West diagonal is a split extension, since it is the pullback of a
split extension.\\
Condition \ref{Farfalle}.iii: To check the commutativity of
$$\xymatrix{E_P \flat H \ar[r]^{\sigma_P  \flat 1} \ar[d]_{1\flat\langle \partial ,p g i \rangle} &
H_0 \flat H \ar[r]^{\xi} & H \ar[d]^{\langle \partial ,p g i \rangle} \\
E_P \flat E_P \ar[rr]_{\chi_{E_P}} & & E_P}$$
compose with the pullback projections $\sigma_P  \colon E_P \to H_0$ and
$\overline{p} \colon E_P \to G_1.$ When composing with $\sigma_P ,$ use
 the naturality of $\chi$ and the precrossed module condition on $\bH.$
When composing with $\overline{p},$ the commutativity of the resulting
diagram easily reduces to condition \ref{Farfalle}.iii on the identity butterfly on $\G.$ \\
Condition \ref{Farfalle}.iv: to check the commutativity of
$$
\xymatrix{
E_P\flat \ar[d]_{1\flat\langle 0,g\rangle}\ar[r]^{\overline{p}d\flat 1}
&G_0\flat G\ar[r]^{\xi}
&G\ar[d]^{\langle 0,g\rangle}
\\
E_P\flat E_P \ar[rr]_{\chi_{E_P}}&&E_P
}
$$
you can either compose once again with the pullback projections, or show that the induced action
$\overline{p}d\flat 1\cdot \xi$ is the unique that makes the subobject $(G,\langle 0,g\rangle)$
closed w.r.t.\ conjugation in $E_P$.
\end{proof}
\begin{SplitFarfalleToMor}\label{SplitFarfalleToMor}{\em
We have just seen in \ref{SplitFarfalle} that  every morphism
$P \colon \bH \to \G$ yields a split butterfly, namely $E_P$. Also the converse is true.
}\end{SplitFarfalleToMor}
Indeed, let
$$\xymatrix{H \ar[rd]^{\kappa} \ar[dd]_{\partial} & &
G \ar[dd]^{\partial} \ar[ld]_{\iota} \\
& E \ar@<0,5ex>[ld]^{\sigma} \ar[rd]_{\rho} \\
H_0 \ar@<0,5ex>[ru]^{s} & & G_0}$$
be a split butterfly. Precomposing the commutative diagram of
condition \ref{Farfalle}.iii with $s \flat 1 \colon H_0 \flat H \to E \flat H$
we get the commutativity of
$$\xymatrix{H_0 \flat H \ar[r]^{\xi} \ar[d]_{s \flat \kappa}
& H \ar[d]^{\kappa} \\
E \flat E \ar[r]_{\chi_{E}} & E}$$
From the universal property of the semi-direct product (see \cite{Janelidze03}, Theorem 1.3), we obtain a unique arrow
$\overline \kappa$ making commutative the diagram
$$\xymatrix{H \ar[r]^{h} \ar[rd]_{\kappa} & H_1 \ar[d]^{\overline \kappa}
& H_0 \ar[l]_{e} \ar[ld]^{s} \\
& E}$$
The requested morphism $P \colon \bH \to \G$ is the one corresponding to the
following internal functor (notation as in \ref{Torsori}, $\Delta$ is the diagonal):
$$\xymatrix{H_1 \ar@<-1ex>[d]_{d} \ar@<1ex>[d]^{c}
\ar[r]^{\overline \kappa} & E \ar@{.>}@<-0,5ex>[d]_{\sigma}
\ar[r]^>>>>>>{\langle 1,\sigma s\rangle} & R[\sigma]
\ar@{.>}@<-0,5ex>[d]_{\sigma_2}
\ar[r]^{\overline \rho} & G_1 \ar@{.>}@<-0,5ex>[d]_{c} \ar[r]^{i} &
G_1 \ar@<-1ex>[d]_{d} \ar@<1ex>[d]^{c} \\
H_0 \ar[u]|{e} \ar[r]_{1} & H_0 \ar@{.>}@<-0,5ex>[u]_{s} \ar[r]_{s} & E
\ar@{.>}@<-0,5ex>[u]_{\Delta} \ar[r]_{\rho} & G_0 \ar@{.>}@<-0,5ex>[u]_{e}
\ar[r]_{1} & G_0 \ar[u]|{e}}$$
(following \cite{CPP}, Proposition 2.1, it suffices to check that this is a morphism of reflexive graphs
and, for this, use the dotted arrows).

\begin{CompRid}\label{CompRid}{\em Given a morphism $Q \colon \K \to \bH$
of  crossed modules and a butterfly $E \colon \bH \to \G,$
we can transform $Q$ into a split butterfly $E_Q \colon \K \to \bH$ as in
\ref{SplitFarfalle} and then to compose
$E_Q$ with $E$ using composition of butterflies described in \ref{CompFarfalle}.
We describe here a somehow easier way to calculate $E_Q \cdot E$ which
will be called reduced composition and denoted by $Q \cdot_{rc} E.$
Starting from
$$\xymatrix{K \ar[r]^{q} \ar[dd]_{\partial} & H \ar[dd]_{\partial}
\ar[rd]^{\kappa} & & G \ar[ld]_{\iota} \ar[dd]^{\partial} \\
& & E \ar[ld]^{\sigma} \ar[rd]_{\rho} \\
K_0 \ar[r]_{q_0} & H_0 & & G_0}$$
consider the pullback
\begin{equation}\label{Diag:reduced_composition}
\raisebox{3ex}{\xymatrix{E' \ar[r]^>>>>{q'} \ar[d]_{\sigma'}
& E \ar[d]^{\sigma} \\
K_0 \ar[r]_{q_0} & H_0}}
\end{equation}
and the arrows
$$\langle 0,\iota\rangle \colon G \to E'
\;\;\;\;\;
\langle \partial, p \cdot \kappa \rangle \colon K \to E'.$$
The  $Q \cdot_{rc} E$ is given by
$$\xymatrix@C=16ex{K \ar[rd]^{\langle \partial, p \cdot \kappa \rangle}
\ar[dd]_{\partial} & &
G \ar[ld]_{\langle 0,\iota\rangle} \ar[dd]^{\partial} \\
& E' \ar[ld]^{\sigma'} \ar[rd]_{q' \cdot \rho} \\
K_0 & & G_0}$$
and it coincides with the butterfly $E_Q\cdot_{rc}E$.
In particular, if $I_{\bH} \colon \bH \to \bH$ is the identity butterfly (\ref{IdFarfalle}),
then $Q \cdot_{rc} I_{\bH}$ is precisely  the split butterfly $E_Q$ as in \ref{SplitFarfalle}.
}\end{CompRid}

\begin{proof} We have to prove that $E_Q \cdot E = Q \cdot_{rc} E$.
Let us consider the following picture, where all the squares are pullbacks and moreover,
down-right square is the discrete fibration  of \ref{RemProfunctors}:
$$
\xymatrix@!=5ex{
&&K_0\times_{q,c} H_1 \times_{d,\sigma}E \ar[dl]_{r}\ar[dr]^{\phi}\\
&E_Q\ar[dl]_{\sigma_Q}\ar[dr]^{\overline{q}}
&&R\ar[dl]_{\overline{\sigma}}\ar@<-.5ex>[dr]_{\overline{d}}\ar@{-->}@<+.5ex>[dr]^{\overline{c}}\\
K_0\ar[dr]_{q_0}&&H_1\ar[dl]^{c}\ar@<-.5ex>[dr]_{d}\ar@{-->}@<+.5ex>[dr]^{c}&&E\ar[dl]^{\sigma}\\
&H_0&&H_0
}
$$
By commutativity of limits, the topmost object is the limit over the $w$-shaped
 diagram $\{ q_0,c,d,\sigma\}$, whence the notation adopted.
The pullback (\ref{Diag:reduced_composition}) determines a unique
$\omega\colon K_0\times_{q,c} H_1 \times_{d,\sigma}E\to E'=K_0\times_{q,\sigma}E$,
such that $\omega {q'}=\phi\,\overline{c}$ and $\omega {\sigma'} = r\sigma_Q$.
Now we can consider the diagram
$$
\xymatrix@!C=16ex{
G\ar[d]_{\langle 0,\iota\rangle}\ar@{}[dr]|{(i)}
& G\ar[d]^{\langle 0,0,\iota\rangle}\ar[l]_1
\\
\ar@{}[dr]|{(ii)}
E' \ar[d]_{{\sigma'}}
&K_0\times_{q,c} H_1 \times_{d,\sigma}E\ar[l]_(.55){\omega}\ar[d]^{r}\ar@{}[dr]|{(iii)}
&H\ar[l]_(.3){\langle0,h,\kappa\rangle}\ar[d]^1
\\
K_0
&E_Q\ar[l]^{\sigma_Q}
&H\ar[l]^(.4){\langle0,h\rangle}
}
$$
By composing with pullback projections, one easily shows that $(i)$ and $(iii)$ commute, so that all the squares
are commutative.  Then, since $r$ is a regular epimorphism, by the (normalized) Barr-Kock Theorem \ref{nBarrKock}, $(ii)$
is a pullback square, hence $\omega$ is a regular epimorphism and it has the same kernel as $\sigma_Q$.
Moreover, since $\ker(\sigma_Q)=\langle0,h\rangle$, $(iii)$ proves that $\ker(\omega)=\langle0,h,\kappa\rangle$

So far, we proved a technical
\begin{Lemma}\label{LemmaOmega}
The sequence
$\xymatrix@C=8ex{H\ar[r]^(.30){\langle0,h,\kappa\rangle}&K_0\times_{q,c} H_1 \times_{d,\sigma}E\ar[r]^(.7){\omega}&E'}$
is exact.
\end{Lemma}
Now we can finally prove reduced composition. To this end, let us consider the following diagram
$$
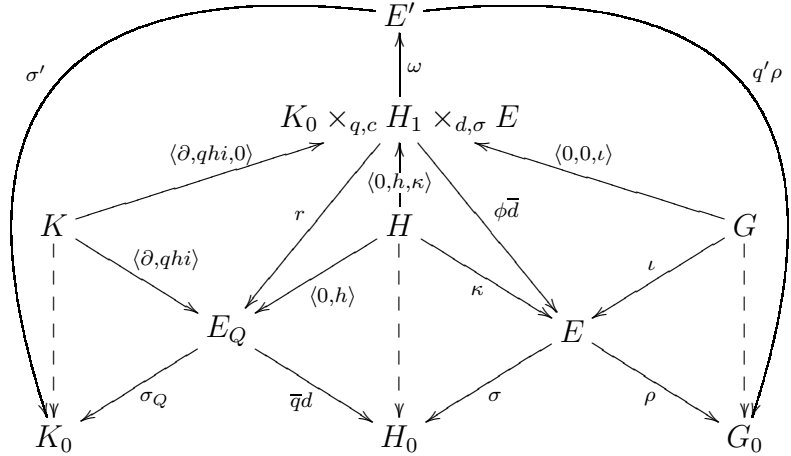
\xymatrix@!C=8ex{
& & E' \ar@/^7pc/[rrdddd]^{{q'}\rho} \ar@/_7pc/[lldddd]_{{\sigma'}}
\\
& & K_0\times_{q,c} H_1 \times_{d,\sigma}E \ar[u]_{\omega} \ar[ldd]_r \ar[rdd]^{\phi \overline{d}}
\\
K \ar[rru]^{\langle \partial,q h i,0 \rangle} \ar@{-->}[dd] \ar[rd]^{\langle \partial,q h i\rangle} & &
H \ar[u]|>>>>>>{\langle 0,h,\kappa \rangle}
\ar[ld]^{\langle 0, h \rangle} \ar[rd]_{\kappa} \ar@{-->}[dd] & &
G \ar[llu]_{\langle 0,0,\iota \rangle} \ar[ld]_{\iota} \ar@{-->}[dd] \\
& E_Q \ar[ld]^{\sigma_Q} \ar[rd]_{\overline{q}  d} & & E \ar[ld]^{\sigma} \ar[rd]_{\rho} \\
K_0 & & H_0 & & G_0}$$
The two butterflies involved are  (from left to right), the split butterfly \mbox{$E_Q\colon\K\to\bH$} corresponding to the
morphism $Q$, and  \mbox{$E\colon \bH\to \G$}.
What we are to show is that the above diagram yields the composition of the two. In fact the resulting butterfly
would be precisely $Q\cdot_{rc}E$, as desired.

By composition of pullbacks, the square $r\cdot\overline{q} d =\phi \overline{d}\cdot \sigma$ above
is a pullback, and by  Lemma \ref{LemmaOmega} $\omega$ is the cokernel of $\langle0,h,\kappa\rangle$.
Moreover $\sigma'$ is (the only morphism) such that $\omega \sigma' = r\sigma_Q$ and ${q'}\rho$ is
(the only one) such that $\omega {q'}\rho=\phi\overline{d}\rho$, and this concludes the proof.
\end{proof}
The following statement will help us in defining the embedding of crossed modules into butterflies.
\begin{Proposition}\label{prop:azione_id_butt}
Reduced composition gives yields on hom-categories a left monoidal action of crossed module morphisms on butterflies, i.e.\ the following formul{\ae} coherently hold when well defined, for morphisms $P$ and $Q$ and for butterflies $E$ and $F$:
\begin{itemize}
\item[$A1.$] $Q\cdot_{rc}EF\cong(Q\cdot_{rc}E)F$
\item[$A2.$] $PQ\cdot_{rc}E\cong P\cdot_{rc}(Q\cdot_{rc}E)$
\item[$A3.$] $I\cdot_{rc}E\cong E$
\end{itemize}
\end{Proposition}
\begin{proof}(outline)
The proof of $A3$ is trivial, that of $A2$ is straightforward. The proof of $A1$ can be easily deduced from
the particular case
\begin{itemize}
\item[$A1^*\!\!.\ $] $Q\cdot_{rc}F\cong (Q\cdot_{rc}I)F$
\end{itemize}
where $I$ is the identity butterfly on the domain of F. Actually  one computes
$$Q\cdot_{rc}EF\cong (Q\cdot_{rc}I)EF \cong ((Q\cdot_{rc}I)E)F \cong ((Q\cdot_{rc}(IE))F \cong (Q\cdot_{rc}E)F$$
Hence we are to prove $A1^*$ holds, but since $Q\cdot_{rc}I=E_Q$, this is precisely the content of the proof
of the consistency of reduced composition described above.
\end{proof}

\section{Butterflies are fractions}

In this section we prove the main result of the paper, but first it is necessary to introduce the fractions whose
the title  refers to. As for the case of groups (see \cite{Noohi05}), given a butterfly it is possible to
construct a span of  morphisms, one  being a weak equivalence.
By denormalizing, this yields a fraction of internal functors.\\

Categories of fractions have been introduced by P. Gabriel and M. Zisman in
\cite{GZ67} to give a simplicial construction of the homotopy category
of CW complexes. In order to study toposes locally equivalent to toposes of
sheaves on a topological space, in \cite{Pronk96} D. Pronk generalized
Gabriel-Zisman concept introducing bicategories of fractions.

\begin{BicatFract}\label{BicatFract}{\em
Imitating the usual universal property of the category of fractions, it is clear
how to state the universal property of the bicategory of fractions
$$\cP_{\Sigma} \colon \cB \to \cB[\Sigma^{-1}]$$
of a bicategory $\cB$ with respect to a class $\Sigma$ of 1-cells (\cite{Pronk96}):
the bicategory of fractions
of $\cB$ with respect to $\Sigma$ is a homomorphism of bicategories
$$\cP_{\Sigma} \colon \cB \to \cB[\Sigma^{-1}]$$
universal among all homomorphisms $\cF \colon \cB \to \cA$ such that $\cF(S)$
is an equivalence for all $S \in \Sigma.$
This means that, for every bicategory $\cA,$
$$\cP_{\Sigma} \cdot - \;\colon \Hom(\cB[\Sigma^{-1}],\cA) \to
\Hom_{\Sigma}(\cB,\cA)$$
is a biequivalence of bicategories, where a homomorphism $\cF \colon \cB \to \cA$
lies in $\Hom_{\Sigma}(\cB,\cA)$ when $\cF(S)$ is an equivalence for all $S \in \Sigma.$
}\end{BicatFract}

The real challenge with bicategories of fractions is to find an explicit, maniable
description of $\cB[\Sigma^{-1}].$ A first general result in this direction, established in
\cite{Pronk96}, states that if $\Sigma$ satisfies some suitable conditions (has a ``right
calculus of fractions'') then the bicategory of fraction exists and can be described as follows:
 the objects of $\cB[\Sigma^{-1}]$ are those of $\cB$
and the 1-cells of $\cB[\Sigma^{-1}]$ are spans of 1-cells in $\cB$ with the backward leg
in $\Sigma$ (this is the non straightforward generalization of a well-known result
from \cite{GZ67}).

In order to prove that butterflies provide the bicategory of fractions of $\GrpdC$ with
respect to weak equivalences, we will use the following result.

\begin{Proposition}\label{PropRCF}(Pronk \cite{Pronk96})
Let $\Sigma$ be a class of 1-cells in a bicategory $\cB$. Assume that $\Sigma$ has a right calculus of
fractions and consider a homomorphism of bicategories $\cF \colon \cB \to \cA$ such that
\begin{enumerate}
\item[EF0.] $\cF(S)$ is an equivalence for all $S \in \Sigma;$
\item[EF1.] $\cF$ is surjective up to equivalence on objects;
\item[EF2.] $\cF$ is full and faithful on 2-cells;
\item[EF3.] For every 1-cell $F$ in $\cA$ there exist 1-cells $G$ and $W$ in $\cB$
with $W$ in $\Sigma$ and a 2-cell $\cF (G) \Rightarrow \cF(W) \cdot F.$
\end{enumerate}
Then the (essentially unique) extension
$$\widehat \cF \colon \cB[\Sigma^{-1}] \to \cA$$
of $\cF$ through $\cP_{\Sigma}$ is a biequivalence.
\end{Proposition}

\begin{CoopXMods}\label{CoopXMods}{\em

Before we show how a butterfly turns into a fraction,
we need one more property of butterflies.
\par
Consider a butterfly $E \colon \bH \to \G$.
The arrows $$\kappa \colon H \to E \leftarrow G \colon \iota$$
{\em cooperate} (see section \ref{HuqIsSmith}), that is, there exists a unique arrow $\varphi=\kappa \sharp \iota$ such that the diagram
$$\xymatrix{H \ar[r]^<<<<<{\langle 1,0\rangle} \ar[rd]_{\kappa}
& H \times G \ar[d]_{\varphi} &
G \ar[l]_<<<<<{\langle 0,1\rangle} \ar[ld]^{\iota} \\
& E}$$
commutes. Indeed, the fact that $\kappa$ and $\iota$ cooperate is equivalent
to the fact that the composition
$$\xymatrix{G \diamond H \ar[r]^{\delta}
& G+H \ar[r]^>>>>>>{[\iota,\kappa]} & E}$$
is the zero morphism (see for example \cite{MM10}), where $\delta$ is the diagonal
of the pullback
$$\xymatrix{G \diamond H \ar[d]_{\delta_2} \ar[r]^{\delta_1}
& G \flat H  \ar[d]^{j_{G,H}} \\
H \flat G \ar[r]_{j_{H,G}} & G+H}$$
The equation $\delta \cdot [\iota,\kappa] = 0$ follows from the commutativity of
$$\xymatrix{G \diamond H \ar[ddd]_{\delta} \ar[r]^{\delta_2} &
H \flat G \ar[r]^{j_{H,G}} & G+H \ar[r]^{[0,1]} & H \simeq 0 \flat H
\ar[d]^{! \flat 1} \\
& & & H_0 \flat H \ar[d]^{\xi} \\
& & & H \ar[d]^{\kappa} \\
G+H \ar[rrr]_{[\iota,\kappa]} & & & E}$$
which can be reduced to the commutativity of
$$\xymatrix{E \flat H \ar[d]_{j_{E,H}} \ar[r]^{\sigma \flat 1} &
H_0 \flat H \ar[r]^{\xi} & H \ar[d]^{\kappa} \\
E+H \ar[rr]_{[1,\kappa]} & & E}$$
Finally, this is a consequence of condition \ref{Farfalle}.iii using that
$$\chi_E = j_{E,E} \cdot [1,1] \colon E \flat E \to E + E \to E$$
}
\begin{Remark} Observe that the fact that $\kappa$ and $\iota$  cooperate may
be used as a starting point for creating many non-trivial examples of butterfly:
one starts by considering two cooperating normal subobjects and then computes their respective cokernels.
\end{Remark}

\end{CoopXMods}
\begin{SpanButterfly}\label{SpanButterfly}
{\em So far we established that the two crossed modules $\kappa$ and $\iota$
cooperate. Now, we are going to prove that $\varphi$ is itself a crossed module,
for a suitable action $\overline{\xi}$, and that the diagram
$$\xymatrix{H \ar[d]_{\partial} & H \times G \ar[l]_>>>>>{\pi_H}
\ar[d]_{\varphi} \ar[r]^>>>>>{\pi_G} & G \ar[d]^{\partial} \\
H_0 \ar@{}[ru]|{(1)} & E \ar[l]^{\sigma} \ar[r]_{\rho}
\ar@{}[ru]|{(2)} & G_0}$$
is  a span of  crossed modules,
$$\xymatrix{\bH & & [E] \ar[ll]_{(\pi_H, \sigma)} \ar[rr]^{(\pi_G, \rho)} & & \G}$$
with  $(1)$ being a weak equivalence.
}\end{SpanButterfly}
\begin{proof}
The commutativity of $(1)$ and $(2)$ can be proved by precomposing with the
jointly epimorphic pair
$$\langle 1,0 \rangle \colon H \to H \times G
\leftarrow G \colon \langle 0,1 \rangle$$
Moreover, $(1)$ is a pullback because it is commutative and the
regular epimorphisms $\pi_H$ and $\sigma$ have same kernel $G$ (use \ref{nBarrKock}).
Therefore, we can apply Lemma \ref{LemmaPBXmod} to $(1)$ getting  that $\varphi$ is a crossed module
 and   that $(1)$ is a weak equivalence of crossed modules.  The action $\overline\xi$ that makes
 $\varphi$ a crossed module is the unique morphism
such that $\overline\xi \pi_H = (\sigma \flat \pi_H)  \xi$
and $\overline\xi  \varphi = (1 \flat \varphi)  \chi_E)$, (see section \ref{LemmaPBXmod}). \\
It remains to show that $(2)$ is a morphism of crossed modules, i.e.\ that the diagram
$$\xymatrix{E \flat (H \times G) \ar[d]_{\overline\xi}
\ar[rr]^{\rho \flat \pi_G} & & G_0 \flat G \ar[d]^{\xi} \\
H \ar@{}[rru]|{(3)} \times G \ar[rr]_{\pi_G} & & G}$$
commutes. For this, we need a different description of $\overline\xi.$\\
Observe that the pullback $(1)$ can be expressed as composite of two
pullbacks using the discrete fibration associated to the left wing of the butterfly
$E \colon \bH \to \G$ (see \ref{Torsori}):
$$\xymatrix{H \times G \ar[r]^>>>>>>{\pi_H} \ar[d]_{\overline h}
\ar@/_2pc/[dd]_{\varphi} & H \ar[d]^{h} \ar@/^2pc/[dd]^{\partial} \\
R \ar[r]^{\overline\sigma} \ar[d]_{d} & H_1 \ar[d]^{d} \\
E \ar[r]_{\sigma} & H_0}$$
Moreover, since $h$ is the kernel of $c \colon H_1 \to H_0$
and $\iota$ is the kernel of $\sigma \colon E \to H_0,$ comparing the diagrams
$$\raisebox{+7ex}{\xymatrix{H \times G \ar[r]^>>>>>>{\pi_H} \ar[d]_{\overline h}
& H \ar[d]^{h} \\
R \ar[r]^{\overline\sigma} \ar[d]_{c} & H_1 \ar[d]^{c} \\
E \ar[r]_{\sigma} & H_0}}
\;\;\; \mbox{ and } \;\;\;
\raisebox{+7ex}{\xymatrix{H \times G \ar[r]^>>>>>>{\pi_H} \ar[d]_{\pi_G} & H \ar[d] \\
G \ar[d]_{\iota} \ar[r] & 0 \ar[d] \\
E \ar[r]_{\sigma} & H_0}}$$
(where all squares are pullbacks) we get the commutativity of
$$\xymatrix{H \times G \ar[r]^>>>>>>{\pi_G} \ar[d]_{\overline h}
& G \ar[d]^{\iota} \\
R \ar[r]_{c} & E}$$
Observe now that since $h \colon H \to H_1$ is a normal mono, there exists a unique
$\chi_{\mid} \colon H_1 \flat H \to H$ such that
$$\xymatrix{
H_1 \flat H \ar[d]_{1 \flat h} \ar[r]^{\chi_{\mid}}
& H \ar[d]^h
\\
H_1 \flat H_1 \ar[r]_{\chi}
& H_1}$$
commutes. From this fact, it follows easily that also
$$\xymatrix{R \flat (H \times G)
\ar[rr]^>>>>>>>>>>>>{(\overline \sigma \flat \pi_H)  \chi_{\mid}}
\ar[d]_{(1 \flat \overline h)  \chi} & & H \ar[d]^{h} \\
R \ar[rr]_{\overline\sigma} & & H_1}$$
commutes. By the universal property of the pullback of $\overline\sigma$ and $h$
we get a unique morphism $x$ such that
$$\xymatrix{R \flat R \ar[d]_{\chi} & & R \flat (H \times G) \ar[d]_{x}
\ar[ll]_{1 \flat \overline h} \ar[rr]^{\overline\sigma \flat \pi_H}
& & H_1 \flat H \ar[d]^{\chi_{\mid}} \\
R & & H \times G \ar[ll]^{\overline h} \ar[rr]_{\pi_H} & & H}$$
commutes. The action $\overline\xi$ factorizes through $x$ as follows:
$$\xymatrix{E \flat (H \times G) \ar[rr]^{\overline\xi} \ar[rd]_{e \flat 1}
& & H \times G \\
& R \flat (H \times G) \ar[ru]_{x} }$$
To check the commutativity of the previous triangle, compose with the pullback
projections
$$\xymatrix{E & \ar[l]_>>>>>>{\varphi} H \times G \ar[r]^>>>>>>{\pi_H} & H}$$
When composing with $\varphi,$ use that $\varphi = \overline{h}  d$
and the left-hand square in the definition of $x.$ When composing with $\pi_H,$
use the right-hand square in the definition of $x$ and the commutativity of
$$\xymatrix{H_0 \flat H \ar[rr]^{\xi} \ar[rd]_{e \flat 1} & & H \\
& H_1 \flat H \ar[ru]_{\chi_{\mid}} }$$
(this last equation is easy to verify: compose with the monomorphism $h$ and use
the definition of $\chi_{H_1}$, see diagram (\ref{Diag:gpd_to_act}) of Proposition \ref{PropGrpdXMod}).\\
We are ready to prove the commutativity of diagram $(3) \colon$replace
$\overline \xi$ by $(e \flat 1)  x,$ compose with the monomorphism
$\iota \colon G \to E$ and use the left-hand square in the definition of $x,$
condition \ref{Farfalle}.iv and the equation $\pi_G  \iota = \overline h  c$
established above.
\end{proof}

\begin{HomF}\label{HomF}

{\em Combining the equivalence
$$J \colon \GrpdC \to \XModC$$
of Proposition \ref{PropGrpdXMod} with the construction of the split butterfly
$E_P$ associated with a morphism $P$
(\ref{SplitFarfalle}), we are ready to define a homomorphism of bicategories
$$\cF \colon \GrpdC \to \cB(\cC).$$

On objects and on 1-cells we let
$$
\cF(\bH)=J(\bH),
\qquad\cF(P \colon \bH \to \G) = (E_{J(P)} \colon J(\bH) \to J(\G)).$$
The composition and the identity structural isomorphisms are defined by using the properties described
in Proposition \ref{prop:azione_id_butt}, by identifying the behavior of $\cF$ on 1-cells with the action of the (reduced) composition with the identity butterfly (see Remark \ref{rem:action_id}).

It remains to define $\cF$ on 2-cells. Let
$\alpha \colon P \Rightarrow Q \colon \bH \to \G$ be a  natural
transformation; there exists a unique morphism $\overline \alpha$ such that the diagram
$$\xymatrix{E_{J(P)} \ar[r]^{\sigma_P} \ar[rd]^{\overline \alpha}
\ar@/_1pc/[rdd]_{\overline{p}} & H_0 \ar[rd]^{\alpha} \\
& G_1 \times_{c,d} G_1 \ar[r]^{\pi_2} \ar[d]_{\pi_1} & G_1 \ar[d]^d \\
& G_1 \ar[r]_c & G_0}$$
commutes. Finally, $\cF(\alpha) \colon E_{J(P)} \to E_{J(Q)}$
is the unique morphism such that the diagram
$$\xymatrix{
E_{J(P)} \ar[r]^<<<<{\overline \alpha} \ar[rd]^{\cF(\alpha)}
\ar@/_1pc/[rdd]_{\sigma_P} & G_1 \times_{c,d} G_1 \ar[rd]^m \\
& E_{J(Q)} \ar[r]^{\overline{q}} \ar[d]_{\sigma_Q} & G_1 \ar[d]^c \\
& H_0 \ar[r]_{Q_0} & G_0}$$
commutes.
Set theoretically, the map $\cF(\alpha)$ sends the pair $(y,\xymatrix{x\ar[r]^(.4)f&p_0(y)})\in E_P$
to the pair $(y,\xymatrix{x\ar[r]^(.4)f&p_0(y)\ar[r]^{\alpha(y)}&q_0(y)})\in E_Q$

}\end{HomF}
\begin{Remark}\label{rem:action_id}{\em
Equivalently, $\cF \colon \GrpdC \to \cB(\cC)$ can be obtained as the composite of $J\colon \GrpdC \to XMod(\cC)$ with
the embedding $\cB\colon XMod(\cC)\to \cB(\cC)$ which is the identity on objects and acts on hom-categories
by the reduced composition with the identity butterfly
$-\cdot_{rc} I_{\G}\colon XMod(\cC)(\bH,\G)\to \cB(\cC)(\bH,\G)$.}
\end{Remark}
\begin{Theorem}\label{ThMain} The homomorphism
$$\cF \colon \GrpdC \to \cB(\cC)$$
defined in \ref{HomF} is the bicategory of fractions of $\GrpdC$ with respect
to the class $\Sigma$ of weak equivalences.
\end{Theorem}

\begin{proof}
Since the class $\Sigma$ has a right calculus of fractions (Propositions 5.5
and 5.2 in \cite{Vitale10}), we have to prove that $\cF$ satisfies conditions EF0 -- EF3
of Proposition \ref{PropRCF}.
\par
EF0 : Consider a weak equivalence of  groupoids and the corresponding
morphism $P \colon \bH \to \G$ of  crossed modules:
$$\xymatrix{H \ar[d]_{\partial} \ar[r]^{p} & G \ar[d]^{\partial} \\
H_0 \ar[r]_{p_0} & G_0}$$
As recalled in Remark \ref{RemWeakEquiv} the arrows induced on kernels and cokernels of $\partial$ are
isomorphisms.
As a first step, we show that the previous diagram
is a pullback. For this, consider the regular epi - mono factorizations of $\partial \colon$
$$\xymatrix{H \ar[r]^{p} \ar[d]_{\partial_1} & G \ar[d]^{\partial_1} \\
I(H) \ar[r]^{I(p)} \ar[d]_{\partial_2}
& I(G) \ar[d]^{\partial_2} \\
H_0 \ar[r]_{p_0} & G_0}$$
By the (normalized) Barr-Kock Theorem \ref{nBarrKock}, the upper square is a pullback because the two regular epimorphisms $\partial_1$
have isomorphic kernels. As far as the lower square is concerned, observe that
$\partial_2 \colon I(H) \to H_0$ is normal (precrossed module condition
in \ref{XMod}) and, therefore, it is the kernel of its cokernel. Using this fact and the fact
that the arrow between cokernels is a monomorphism, it is easy to check that the lower
square satisfies the universal property of the pullback. \\
Now, we want to show that the split butterfly
$$E_P \colon \bH \to \G$$
associated to the above morphism of  crossed modules as in
\ref{SplitFarfalle} is an equivalence. Following Proposition \ref{PropEquivFarfalle},
it is enough to show that  $E_P$ is flippable.
For this, consider the diagram
$$\xymatrix{H \ar[d]_{h^{\bullet}} \ar[r]^{p} & G \ar[dd]^{g^{\bullet}} \\
H_1 \ar[d]_{\langle d, p_1 i \rangle} \\
E_P \ar[d]_{\sigma_P } \ar[r]^{\overline{p}} & G_1 \ar[d]^{c} \\
H_0 \ar[r]_{p_0} & G_0}$$
The whole square is precisely $P \colon \bH \to \G,$ so that it is a pullback.
The lower square also is a pullback (\ref{SplitFarfalle}), so that the upper square is a
pullback. From this and the fact that $g^{\bullet}$ is the kernel of
$d \colon G_1 \to G_0,$ we immediately get that $h^{\bullet}  \langle d,p_1 i \rangle$
is the kernel of $\overline{p}  d \colon E_P \to G_0.$
Finally, $\overline{p} d$ is a regular epimorphism by definition of essential
surjective (\ref{WealEquiv}) and, therefore, it is the cokernel of its kernel.
\par
EF1: Since $\cF$ on objects is the composite
$$\GrpdC \to \XModC \to \cB(\cC)$$
with the first step being an equivalence and the second one being the identity on objects,
condition EF1 is clearly satisfied.
\par
EF2: We are to prove that $\cF \colon \GrpdC \to \cB(\cC)$ is full and faithful
on 2-cells. To this end, let us consider two parallel morphisms of  crossed modules
$$P, Q \colon \bH \rightrightarrows \G$$
and a morphism
$$f \colon E_P \to E_Q$$
between the corresponding split butterflies (\ref{SplitFarfalle}), i.e.\ the following four triangles commute:
$$
\xymatrix@=7ex{\ar@{}[dr]|(.25){(i)}|(.75){(ii)}
H\ar[r]^{\langle\partial,p g^{\bullet}\rangle}\ar[d]_{\langle\partial,q g^{\bullet}\rangle}
&E_P\ar[dl]|f \ar[d]^{\overline{p}d}
\\
E_Q\ar[r]_{\overline{q}d}
&G_0
}
\qquad
\xymatrix@=7ex{
E_P\ar[d]_{\sigma_P}\ar[dr]|f
&G\ar[l]_{\langle0,g\rangle}\ar[d]^{\langle0,g\rangle}
\ar@{}[dl]|(.25){(iv)}|(.75){(iii)}\\
H_0&E_Q\ar[l]^{\sigma_Q}
}$$
Consider also the arrow $\pi$ given by the universal property of the pullback $E_P \colon$
$$\xymatrix{H_0 \ar[r]^{e} \ar@/_1,5pc/[rdd]_{1} \ar[rd]^{\pi}
& H_1 \ar@/^1pc/[rd]^{p_1} \\
& E_P \ar[d]_{\sigma_P } \ar[r]^{\overline{p}} & G_1 \ar[d]^{c} \\
& H_0 \ar[r]_{p_0} & G_0}$$
Put
$$\alpha_f \colon \xymatrix{H_0 \ar[r]^{\pi} & E_P \ar[r]^{f} &
E_Q \ar[r]^{\overline{q}} & G_1}$$
It is easy to check that $\alpha_f \cdot d = p_0$ and $\alpha_f \cdot c = q_0$: just use commutativity of
$(ii)$ and $(iii)$ above. To prove that $\alpha$ is natural requires some computations.
Following the characterization of Proposition \ref{PropNatPFFGraph}, we are to show that
$$
\xymatrix{
H\ar[r]^{\partial}\ar[d]_{\langle p,q\rangle}
&H_0\ar[r]^{\pi}
&E_P\ar[r]^{f}
&E_Q\ar[d]^{\overline q}
\\
G\times G\ar[rrr]_{m_0}
&&&G_1}
$$
commutes, where $m_0=g\sharp g^{\bullet}$ is the cooperator of $g$ and $g^{\bullet}$.
To this end, let us consider the following diagram, whose outer rectangle gives naturality of $\alpha_f$:
$$
\xymatrix@C=10ex{
H\ar[r]^{\partial\pi}\ar[dr]_{\langle1,p\rangle}\ar[dd]_{\langle q,p\rangle}
&E_P\ar[r]^f\ar@{}[dl]|(.3){(1)}\ar@{}[dr]|(.3){(2)}
\ar@{}[ddl]\ar@{}[ddr]|(.5){(3)}\ar@{}[dd]|(.8){(4)}
&E_Q\ar[dd]^{\overline q}
\\
&H\times G \ar[u]_{\varphi_P}\ar[ur]_{\varphi_Q}\ar[dr]|{qg^\bullet\sharp\, g}\ar[dl]^{\langle q,1\rangle}&&
\\
G\times G \ar[r]_{\textrm{twist}}
&G\times G\ar[r]_{m_0}
&G_1}
$$
where the maps $\varphi_P$ and $\varphi_Q$ are the cooperators relative to butterflies $E_P$ and $E_Q$ (see \ref{CoopXMods}).
The commutativity of $(2)$, $(3)$ and $(4)$ is easily obtained by uniqueness of cooperators, by means of the precompositions
with canonical morphisms $\xymatrix{H\ar[r]^(.4){\langle1,0\rangle}&H\times G&G\ar[l]_(.4){\langle0,1\rangle}}$. Observe that
in proving $(2)$ we use precisely the hypothesis $(i)$ and $(iv)$ above. In fact $(1,f)$ is precisely
the morphism of the spans (determined by the butterflies $E_P$ and $E_Q$) corresponding to $f$.
Finally we show that $(1)$ commutes, by composing with pullback projections $\sigma_P$ and $\overline p$:
\begin{itemize}
\item[1.]
$\partial\cdot\pi\cdot \sigma_P=\partial=\langle 1,p\rangle \cdot\mathrm{pr}_1\cdot\partial=\langle 1,p\rangle\cdot\varphi_P\cdot\sigma_P$
where the last equality is just the weak equivalence $(\mathrm{pr}_1,\sigma_P)\colon \varphi_P\to \partial$ in the span
of crossed modules corresponding to $E_P$.
\item[2.] First observe that
$$
\xymatrix{
G\ar[rr]^{\partial e}\ar[dr]_{\langle1,1\rangle}&&G_1\\
&G\times G\ar[ur]_{m_0}
}
$$
commutes. This can be easily deduced from the very definition of $m_0$  (this equation
is one of the axioms defining a Peiffer graph, see \cite{MM10}).
Hence we start our computation: \\
$\partial \cdot\pi\cdot \overline{p} = \partial \cdot p_0\cdot e = p\cdot \partial\cdot e
=\langle p,p\rangle \cdot m_0 =\langle1,p\rangle \cdot (p\times 1)\cdot m_0 =
 \langle1,p\rangle \cdot\varphi_P\cdot\overline{p}$.
The last equality is obtained by observing that both $(p\times 1)\cdot m_0$ and $\varphi_P\cdot\overline{p}$
are {\em the} cooperator of $pg$ and $g^{\bullet}$.
\end{itemize}
\par
EF3: We want to prove that the diagram in $\cB(\cC)$
$$\xymatrix{ & [E] \ar[ld]_{E_{(\pi_H,\sigma)}} \ar[rd]^{E_{(\pi_G,\rho)}} \\
\bH \ar[rr]_{E} & & \G}$$ commutes (up to 2-cell).
For this, we use reduced composition described in \ref{CompRid}, and we compute
$(\pi_H,\sigma)\cdot_{rc} E=$
$$
\raisebox{+9ex}{\xymatrix{H \times G \ar[rr]^{\pi_H}
\ar[rd]^{\langle\varphi,\pi_H\kappa\rangle} \ar[dd]_{\varphi}
& & H \ar[dd]_<<<<<{\partial} \ar[rd]^{\kappa}
& & G \ar[ld]_{\iota} \ar[dd]^{\partial} \\
& R[\sigma] \ar[ld]_{\sigma_1} \ar[rr]^<<<<<<{\sigma_2}
& & E \ar[ld]^{\sigma} \ar[rd]_{\rho} \\
E \ar[rr]_{\sigma} & & H_0 & & G_0}}=
\raisebox{+9ex}{\xymatrix{H\times G
\ar[rd]^<<<<<<{\langle\varphi,\pi_H\kappa\rangle} \ar[dd]_{\varphi}
\ar[dd]_{\varphi} & & G
\ar[ld]_>>>>>{\langle 0,\iota\rangle} \ar[dd]^{\partial} \\
& R[\sigma] \ar[ld]^{\sigma_1} \ar[rd]_{\sigma_2  \rho} \\
E & & G_0}}$$
and $(\pi_G,\rho)\cdot_{rc} I_{\G} =$
$$
\raisebox{+9ex}{\xymatrix{H \times G \ar[rr]^{\pi_G}
\ar[rd]^{v} \ar[dd]_{\varphi}
& & G \ar[dd]_<<<<<{\partial} \ar[rd]^{g^{\bullet}}
& & G \ar[ld]_{g} \ar[dd]^{\partial} \\
& R[\sigma] \ar[ld]_{\sigma_1} \ar[rr]^<<<<<<{\overline{\rho} i}
& & G_1 \ar[ld]^{c} \ar[rd]_{d} \\
E \ar[rr]_{\rho} & & G_0 & & G_0}}=
\raisebox{+9ex}{
\xymatrix{H \times G \ar[dd]_{\varphi}
\ar[rd]^{v}
& & G \ar[ld]_{u} \ar[dd]^{\partial} \\
& R[\sigma] \ar[ld]^{\sigma_1} \ar[rd]_{\overline{\rho} i  d} \\
E & & G_0}}$$
where the diagram  $\sigma_1 \cdot \rho = \overline{\rho}i \cdot c$ is a pullback as far as the diagram
$\sigma_1\cdot\rho = \overline{\rho}\cdot d$ is (compare with the diagram in \ref{RemProfunctors}).
We are going to prove that the two butterflies obtained this way are isomorphic. In fact they coincide.\\
We see that $\overline{\rho} i  d=\overline{\rho} c = \sigma_2 \rho$ commutes as it is
one of the pullbacks in the diagram of Remark \ref{RemProfunctors}.
To show that $\langle 0,\iota\rangle = u$ it suffices to compose the last with the
pullback projections $\sigma_1$ and $\sigma_2$ that define $\langle 0,\iota\rangle$:
$u \sigma_1=0$ by definition, $u \sigma_2$ is the normalization of the groupoid
$(R[\sigma],\sigma_1,\sigma_2)$, hence the morphism $\iota$.\\
In conclusion we prove that ${\langle\varphi,\pi_H\kappa\rangle}=v$.
For this it suffices to compose with the second projection $\overline{\rho} i$, since the first is identical.
Hence we consider the following commutative diagrams:
$$\xymatrix{H \times G \ar[r]^<<<<<{\pi_H} \ar@{}[rd]|{(4)}
\ar[d]_{\langle \varphi, \pi_H  \kappa \rangle} & H \ar[d]^{\kappa} \\
R[\sigma] \ar[r]^<<<<<<{\sigma_2} \ar[d]_{\sigma_1}
\ar@{}[rd]|{(5)} & E \ar[d]^{\sigma} \\
E \ar[r]_{\sigma} & H_0}
\;\;\;\;\;\;
\xymatrix{H \times G \ar[r]^<<<<<{\pi_H} \ar@{}[rd]|{(4)}
\ar[d]_{\langle \varphi, \pi_H  \kappa \rangle} & H \ar[d]^{\kappa} \\
R[\sigma] \ar[r]^<<<<<<{\sigma_2} \ar[d]_{\overline{\rho}\cdot i}
\ar@{}[rd]|{(6)} & E \ar[d]^{\rho} \\
G_1 \ar[r]_{d} & G_0}$$
Since $\langle \varphi, \pi_H  \kappa \rangle  \sigma_1 = \varphi$
and $\kappa  \sigma = \partial,$ then $(4)+(5)$ is a pullback. Since
$(5)$ is a pullback and $\sigma_2$ is a split epimorphism, $(4)$ is a pullback.
But $(6)$ also is a pullback, so that $(4)+(6)$ is a pullback.
Since $\kappa  \rho = 0$ and $g  i=g^{\bullet}$ is the kernel of $d,$
$(4)+(6)$ can be written as
$$\xymatrix{H \times G \ar[r]^<<<<<{\pi_H} \ar[d]_{\pi_G} & H \ar[d] \\
G \ar[r] \ar[d]_{g^{\bullet}} & 0 \ar[d] \\
G_1 \ar[r]_{d} & G_0}$$
In particular, this gives
 $\langle \varphi, \pi_H  \kappa \rangle\,   \overline{\rho}\, i
= \pi_G\,  g^{\bullet},$ as requested.
\end{proof}

\begin{Remark}\label{RemSplitToButter}{\em
Condition EF3  above gives an alternative proof of the construction of section \ref{SplitFarfalleToMor}.
Moreover it yields a recipe to obtain the crossed modules morphism corresponding to a split butterfly.

Let us consider the span associated to the split butterfly $E=(\kappa,\iota,\sigma,\rho)$, and suppose we have
chosen a  section $s$ of the split epimorphism $\sigma$. Then, since the left leg of the span is a pullback diagram,
we can pull the section $s$ back along $\varphi$, and get the morphism of crossed modules $(\overline{s},s)$, moreover
the last is a section of $(\pi_H,\sigma)$.
The situation is summarized in the diagram below.
$$\xymatrix{\ar@<+1ex>[r]^(.45){\overline{s}}
H \ar[d]_{\partial} & H \times G \ar[l]^(.55){\pi_H}
\ar[d]_{\varphi} \ar[r]^(.55){\pi_G} & G \ar[d]^{\partial} \\
H_0\ar@<+1ex>[r]^s & E \ar[l]^{\sigma} \ar[r]_{\rho}
 & G_0}$$
Then we can compose on the left the equation $E_{(\pi_H,\sigma)}E=E_{(\pi_G,\rho)}$ with the morphism
$(\overline{s},s)$ and get:
\begin{eqnarray*}
E&=&(\overline{s},s)(\pi_H,\sigma)\cdot_{rc}E=(\overline{s},s)\cdot_{rc}E_{(\pi_H,\sigma)}E
=(\overline{s},s)\cdot_{rc}E_{(\pi_G,\rho)}=\\
&=&(\overline{s},s)(\pi_G,\rho)\cdot_{rc}I_{\G}=
(\overline{s}\pi_G,s\rho)\cdot_{rc}I_{\G}=E_{(\overline{s}\pi_G,s\rho)}
\end{eqnarray*}
i.e.\ the morphism $(\overline{s}\pi_G,s\rho)$ is associated with the (split) butterfly $E$. Notice
the arbitrary choice of the section $s$: if another section is chosen, the construction yields another
representant of $E$ in the same 2-isomorphism class.
}\end{Remark}

\section{From butterflies to weak morphisms: three concrete examples}
In the following we show how to construct the weak morphism associated to a butterfly in the cases of groups,
Lie algebras and Rings (but the technique can be adapted to other semi-abelian algebraic varieties), and we give an
idea of how to recover a butterfly from a weak morphism.
 The first two instances  are well known in the literature: a butterfly between  two crossed modules of
groups corresponds to a weak morphism of strict 2-groups (from Theorem \ref{ThMain} and \cite{Vitale10}), similarly
a butterfly in Lie algebras gives a homomorphism (i.e.\ semi-strict morphism) of strict Lie 2-algebras
(see \cite{Noohi05} and \cite{NoohiLie}).
In a similar way, a butterfly in rings provides the data for a weak morphism of strict 2-rings. This notion seems not to be present
in the literature, although  there are two notions of weak 2-rings with units (the categorical rings of \cite{catrings} and
the Ann-categories of \cite{anncat}). These two notions coincide in the strict case,  and it is possible to show that our weak morphism of strict 2-rings
specializes to those.

\subsection{The technique}
Let us consider the butterfly $E=(E,\kappa,\rho,\iota,\sigma)$
in a semi-abelian algebraic variety $\mathcal{C}$, and let $U:\mathcal{C}\to \mathcal{S}$
(the axiom of choice holding in $\mathcal{S}$) a suitable  functor that forgets part of the structure.
Let $s$ be a section of $U(\sigma)$.
$$\xymatrix{H \ar[rd]^{\kappa} \ar[dd]_{\partial} & &
G \ar[dd]^{\partial} \ar[ld]_{\iota} \\
& E \ar@<0,5ex>[ld]^{\sigma} \ar[rd]_{\rho} \\
H_0 \ar@<0,5ex>@{-->}[ru]^{s} & & G_0}$$
We want to show how $E$ yields a weak morphism of  groupoids $F_E\colon \bH \to \G$.
From now on we will write just $\sigma$ for $U(\sigma)$, etc.

The functor $U$ preserves finite limits, so that it extends to a 2-functor between the 2-categories of internal groupoids.
Now, to the butterfly $E$ is associated a span (see section \ref{SpanButterfly}) in $\Grpd(\mathcal{C})$ with the
left leg being a weak equivalence:

$$
\xymatrix{\bH & [E]\ar[l]_S\ar[r]^R&\G},
$$
more explicitly:
$$\xymatrix{H_1 \ar@<-0,5ex>[d]_{d} \ar@<0,5ex>[d]^{c} & &
 (H \times G) \rtimes_{\overline \xi} E
 \ar[ll]_>>>>>>>>>>{\pi_H \rtimes \sigma} \ar[rr]^>>>>>>>>>>{\pi_G \rtimes \rho}
 \ar@<-0,5ex>[d]_{d} \ar@<0,5ex>[d]^{c} & &
 G_1 \ar@<-0,5ex>[d]_{d} \ar@<0,5ex>[d]^{c} \\
 H_0 & & E \ar[ll]^{\sigma} \ar[rr]_{\rho} & & G_0}
$$

By applying $U$ to this construction, $S$ turns in an equivalence  in $\Grpd(\mathcal{S})$, so that
it has a weak inverse $S^*$.

The composition $S^*R$ (which is an internal functor in $\Grpd(\mathcal{S})$) is a good candidate for a
weak morphism in  $\Grpd(\mathcal{C})$, with the coherence conditions encoded in the short exact
sequence of the butterfly.

\subsection{Case study: groups}
Let $\mathcal{C}=\Grp$, and $U:\Grp\to\Set_{*}$ the underlying pointed-set functor.

Under the equivalence between crossed modules and groupoids, the crossed module $\partial\colon H\to H_0$
gives rise to the groupoid in groups
$$
\xymatrix@C=10ex{G_1\ar@<+1ex>[r]^d\ar@<-1ex>[r]_c&G_0\ar[l]|e}
$$
where $G_1$ is the semidirect product $G\rtimes G_0$, and structure maps result (additive notation)
$$
c\colon (g,x)\mapsto x, \qquad d\colon (g,x)\mapsto \partial g + x, \qquad e\colon x\mapsto (0,x).
$$
Define the monoidal functor $F_E=(F_0,F_1,F_2)$:
\begin{eqnarray*}
  F_0 &=&s\rho \colon H_0\to G_0;\qquad x\mapsto \rho(s x)  \\
  F_1&=&F\rtimes F_0 \qquad \mathrm{where} \\
  &&F \colon H \to G;\qquad h\mapsto -\kappa(h)+s(\partial(h))\\
  F_2 &\colon& H_0\times H_0\to G_1; \qquad(x,y) \mapsto (sx +sy-s(x+y),\rho(s(x+y))
\end{eqnarray*}
Notice that, since $\partial(sx +sy-s(x+y)) + \rho(s(x+y) = \rho((sx +sy-s(x+y))) + \rho(s(x+y) = \rho(sx)+\rho(sy)$,
$ F_2(x,y)$  is to be interpreted as an arrow $F_0(x+y)\to F_0(x)+F_0(y)$.

From the classification of group extensions (see [ML Homology], for instance) we know that with the short exact sequence
$$
\xymatrix{G\ar[r]^{\kappa}&E\ar[r]^{\sigma}&H_0}
$$
with a chosen set-theoretical section $s$ of $\sigma$ we can associate two of functions
\mbox{$\alpha\colon H_0\to \mathrm{Aut} G$} and \mbox{$f\colon H_0\times H_0\to G$}:
with $\alpha(x)(g)=x\cdot g =sx+g-sx$ and $f(x,y)=sx +sy-s(x+y)$. Such functions
satisfy the following well known relation: for any $x,y,z$ in $H_0$
$$
x\cdot f(y,z)+f(x,yz)=f(x,y)+f(xy,z).
$$
It is now easy to show that this relation corresponds precisely to what is necessary in order
to prove (associative) coherence for the monoidal functor $F_E$.

A further remark. In giving the example above we stressed the role of group extensions classification.
This is surely of interest, but it is not
conceptually necessary. Actually, the reason why we can extend the method used in Remark \ref{RemSplitToButter},
is that the monoidal functor $S$ above, being a monoidal equivalence, has a weak inverse.

Finally we give a glance to the construction of a butterfly from a (normalized) monoidal functor of 2-groups.
Consider, $F=(F_0,F_1)\colon \bH\to \G$ a functor with monoidal structure isomorphisms $F_2^{x_1,x_2}\colon F_0x_1+F_0x_2\to F_0(x_1+x_2)$.
Define $P_0$  as the following limit (in Set)
 $$\xymatrix{ & & P_0 \ar[lld]_{\sigma} \ar[d] \ar[rrd]^{\rho} \\
 H_0 \ar[rd]_{F_0} & & G_1 \ar[ld]^{c} \ar[rd]_{d}  & & G_0 \ar[ld]^{1} \\
 & G_0 & & G_0}$$
i.e.
$$P_0=\{(y,g,x)\in H_0\times G_1 \times G_0\ \mathrm{s.t.}\ F_0y=cg \ \mathrm{and}\ dg=  x \}$$
and
$$
\sigma\colon (y,g,x)\mapsto y, \qquad \rho\colon (y,g,x)\mapsto x
$$
Despite the fact that $F_0$ is not a group homomorphism, it can be proved that
$P_0$ is a group and $\sigma$ and $\rho$ are group homomorphisms (see \cite{Vitale10}, Proposition 6.3).
Moreover, $\ker \sigma = \ker c$ (just because $c$ is surjective).
 Finally, the factorization $\kappa \colon H \to P_0$ of the commutative diagram
 $$\xymatrix{ & & H \ar[lld]_{hd} \ar[d]_{hiF_1} \ar[rrd]^{0} \\
 H_0 \ar[rd]_{F_0} & & G_1 \ar[ld]^{c} \ar[rd]_{d}  & & G_0 \ar[ld]^{1} \\
 & G_0 & & G_0}$$
 through the limit $P_0$ is also a group homomorphism (this follows immediately from
 the naturality of the monoidal structure of $F).$ We get in this way the required butterfly
 $(P_0,\kappa,\rho,\iota,\sigma)$.

\subsection{Case study: Lie algebras}
A groupoid in $\Lie$ is called a strict Lie 2-Algebra in \cite{BC04}. Now we consider the forgetful functor
$U:\Lie\to\Vect$. As for the case of groups, it extends to internal groupoids.
With notation as above we can define $F_E$ with the same technique. Indeed $F_0$ and $F_1$ are defined in the same
way (provided the semidirect product is performed in $\Lie$!), while
$$
F_2\colon (x,y)\mapsto ([sx,sy]-s[x,y],\rho(s[x,y])).
$$
From the theory of Lie algebras extensions, we know that with the extension $(\iota,\sigma)$ (and a chosen linear
 section $s$ of $\sigma$) is associated a linear map $\alpha\colon H_0\to \Der G$, $\alpha(x)(g)=x\cdot g=[sx,g]$, and
a bilinear skew-symmetric map $f:H_0\times H_0\to G$, $f(x,y)= [sx,sy]-s[x,y]$. These maps satisfy the relations
\begin{itemize}
\item[(i)] for any $x,y$ in $H_0$, $[\alpha(x),\alpha(y)]- \alpha([x,y])=\mathrm{ad}_{f(x,y)}$
\item[(ii)] for any $x,y,z$ in $H_0$
$$\sum_{\mathrm{cyclic}}(x\cdot f(y,z)-f([x,y],z))=0$$
\end{itemize}
where $\mathrm{ad}_g$ is the (adjoint) action defined by  $\mathrm{ad}_g(g')=[g,g']$.
Computations show that the first relation helps in proving the naturality of $F_2$, while
the second relation yields the coherence of the bracket operation with respect to the jacobian identity.

\subsection{Case study: rings}
We call  (strict) \emph{$2$-ring} a  groupoid in the category of rings. There are two obvious
forgetful 2-functors from the 2-categories of (strict) categorical rings \cite{catrings} and of (strict)
Ann-categories \cite{anncat} respectively, both with strict homomorphisms and 2-homomorphisms as 1-cells and 2-cells.

We consider the forgetful functor $U\colon \Rng \to \Set_{*}$, that, as for the previous cases,
extends to groupoids. The definition of $F_E$ goes verbatim as in the case of groups, thanks to the additive notation
used there, that now expresses the underlying abelian group structure of a ring.

Again, the exact sequence $(\iota,\sigma)$  provides the data
for proving that $F_E$ is a 2-ring homomorphism. In fact we use $s$, the set-theoretical section of $\sigma$,
to define $f,\epsilon\colon H_0\times H_0\to G$: $f(x,y)=sx+sy-s(x+y)$,
$\epsilon(x,y)=sx\cdot sy-s(x\cdot y)$, and a map $\alpha\colon H_0\to \mathrm{Bim} G$ with
$\alpha(x)(g)=(sx\cdot g,g\cdot sx)$. Then the following relations hold for any $x,y,z$ and $t$ in $H_0$
 (see \cite{ML58}):
\begin{itemize}
\item[(i)] $\alpha(x)+\alpha(y)-\alpha(x+y)=\mu_{f(x,y)}$

\item[(ii)] $\alpha(x)\circ\alpha(y)-\alpha(xy)=-\mu_{\epsilon(x,y)}$

\item[(iii)] $f(0,y)=0=f(x,0)$ and $\epsilon(0,y)=0=\epsilon(x,0)$

\item[(iv)] $f(x,y)+f(z,t)-f(x+z,y+t)-f(x,z)-f(y,t)+f(x+y,z+t)=0$

\item[(v)] $-\epsilon(x,t)-\epsilon(y,t)+\epsilon(x+y,t)+f(xt,yt)-f(x,y)\cdot t=0$

\item[(vi)] $\epsilon(t,x)+\epsilon(t,y)-\epsilon(t,x+y)-f(tx,ty)+f\cdot h(x,y)=0$

\item[(vii)] $x\cdot \epsilon(y,z)-\epsilon(xy,z)+\epsilon(x,yz)-\epsilon(x,y)\cdot z=0$

\end{itemize}
where $\mu_g$ is the inner bimultiplication induced by the multiplication with $g$.

Now,  (i) and (ii) give the naturality  of $F_2$. Moreover, since the normalization conditions
(iii) hold, the relation $(iv)$ gives at once associative and symmetric coherence: actually for
$y=0$ we obtain the cocycle condition for the underlying (abelian) group extension, while letting $x=t=0$ we get
the symmetric coherence.
Finally (vii) yields the associative coherence for the multiplication, and (v) and (vi) give
the distributive coherence.

\begin{Remark}\label{Anafuntori}{\em
If $\cC$ is a category with finite limits (not necessarily semi-abelian),
Theorem \ref{ThMain} may fail. In this  more general case, internal categories may differ from
internal groupoids and the bicategory of fractions $$\CatC[\Sigma^{-1}]$$ may still admit
a (more involved) explicit description: it is the bicategory of internal anafunctors.
This has been proved independently by M. Dupont
in \cite{Dupont08}, where the base category $\cC$ is assumed to be regular,
and by D. Roberts in \cite{Roberts10}, where essential surjectivity is intended
relatively to a Grothendieck topology on $\cC$ and internal categories (not only
internal groupoids) are considered.
}\end{Remark}

\section{Classification of extensions}

In this section we assume that $\cC$ has split extensions
classifiers (see \cite{Bourn06}, and section \ref{sec:semiab}), as it happens, for instance, in the category of groups or of Lie-algebras.
 Consider two objects $H$ and $G$ in $\cC.$
Let $D(H)=(0 \to H)$ be the discrete crossed module on $H$ and
$$\cA(G)=(\cI_G \colon G \to \Aut G \,,\; \ev \colon \Aut G \flat G \to G)$$
the  crossed module associated with the split extensions classifier $\Aut G$
(that is, the  crossed module corresponding to the action groupoid).\\
The following lemma generalizes Example 13.4 of \cite{Noohi05}.
\begin{Lemma}\label{LemmaSchreier}
The groupoid
$$\Ext(H,G)$$
of extensions of the form $H \leftarrow E \leftarrow G$ is isomorphic to the groupoid
$$\cB(\cC)(D(H),\cA(G))$$
Such an isomorphism restricts to split extensions and split butterflies.
\end{Lemma}

\begin{proof}
Let us start with a butterfly
$$\xymatrix{0 \ar[rd]^{!} \ar[dd]_{!} & & G \ar[ld]_{\iota} \ar[dd]^{\cI_G} \\
& E \ar[ld]^{\sigma} \ar[rd]_{\rho} \\
H & & \Aut G}$$
We are going to prove that $\rho$ is uniquely determined.
Following Remark \ref{RemEquivarianzaFarf}, the right wing determines a discrete fibration of
groupoids. Hence the diagram
$$
\xymatrix{\ar@{}[drr]|{(*)}
G \rtimes_{(\rho \flat 1) \ev}E \ar[rr]^{p_E} \ar[d]_{1 \rtimes \rho}
 & & E \ar[d]^{\rho} \\
G \rtimes_{\ev}\Aut G \ar[rr]_{p_{\Aut G}} & & \Aut G}
$$
is a pullback. Therefore, $\rho$ is the unique arrow making the diagram above a pullback
(universal property of $\Aut G).$ \\
Conversely, consider an extension
$$\xymatrix{H & E \ar[l]_{\sigma} & G \ar[l]_{\iota}}$$
Since $\iota$ is normal in $E,$ there exists a unique $\chi_{\mid}$ such that
$$\xymatrix{E \flat G \ar[r]^{\chi_{\mid}} \ar[d]_{1 \flat \iota}
& G \ar[d]^{\iota} \\
E \flat E \ar[r]_{\chi_E} & E}$$
commutes. By the universal property of $\Aut G,$ we get a unique $\rho$ such that
diagram $(*)$ above is a pullback and
$$\xymatrix{E \flat G \ar[rr]^{\chi_{\mid}} \ar[rd]_{\rho \flat 1}
& & \Aut G \flat G \\
& G \ar[ru]_{\ev} }$$
commutes. It remains to show that the extension $(\iota,\sigma),$ equipped
with $\rho,$ is a butterfly from $D(H)$ to $\cA(G).$ Condition \ref{Farfalle}.iv
follows by pasting together the diagram $(*)$ with the following one:
$$
\xymatrix{\ar@{}[drr]|{(**)}
G\rtimes_{\chi_G}G\ar[d]_{1\rtimes\iota}\ar[rr]^{P_G}& &G\ar[d]^{\iota}
\\
G \rtimes_{(\rho \flat 1) \ev}E \ar[rr]_{p_E}
 & & E }
$$
As far as the commutativity
of the right wing is concerned, observe that diagram $(**)+(*)$  is a pullback
(because both $(*)$ and $(**)$ are pullbacks) and use once again the universal property
of $\Aut G$ in order to conclude that $\iota  \rho = \cI_G.$
\end{proof}

\begin{Classification}\label{ClassExt}{\em
Combining the previous isomorphism of groupoids with Theorem \ref{ThMain},
we get a very general classification of extensions:
$$\Ext(H,G) \simeq \cB(\cC)(D(H),\cA(G)) \simeq
\GrpdC [\Sigma^{-1}](D(H),\cA(G))$$
To recover the classical classification of group extensions in terms of factor sets due to
Schereier it suffices to use another result from \cite{Vitale10} already quoted in the
Introduction: when $\cC$ is the category of groups, the bicategory of fractions
$\GrpdC [\Sigma^{-1}]$ can be described as the bicategory of  groupoids,
monoidal functors and monoidal natural transformations, and a monoidal functor from
$D(H)$ to $\cA(G)$ is nothing but a factor set.

}\end{Classification}

\section{The free exact case}

When $\cC$ is the category of groups, the main result of \cite{Noohi05}
is not stated in terms of bicategory of fractions, but it is stated as an equivalence
of groupoids
$$\cB(\cC)(\bH,\G) \simeq \XModC(\bK,\G)$$
where $\bK$ is the crossed module of groups  obtained from $\bH$ by pulling back
$\partial \colon H \to H_0$ along a surjective homomorphism $K_0 \to H_0$
with $K_0$ being a free group. The same is done for Lie algebras in \cite{Abbad10}.
The aim of this section is to generalize the previous equivalence to the case when the
base category $\cC$ is also free exact.

\begin{FreeEx}\label{FreeEx}{\em
We assume that the semi-abelian category $\cC$ is free exact in the sense of
\cite{CV98}, that is, it has enough regular projective objects. This means that for
every object $X$ in $\cC$ there exists a regular epimorphism $x \colon X' \to X$
with $X'$ regular projective. All semi-abelian varieties are of this kind.
In particular, groups and Lie algebras are free exact semi-abelian categories.
}\end{FreeEx}

\begin{Replacement}\label{Replacement}{\em
Let $\C$ be a  groupoid and $s_0 \colon X_0 \to C_0$ a regular
epimorphism with $X_0$ regular projective. Consider the limit
$$\xymatrix{ & & X_1 \ar[lld]_{d} \ar[d]^{s_1} \ar[rrd]^{c} \\
X_0 \ar[rd]_{s_0} & & C_1 \ar[ld]^{d} \ar[rd]_{c} & & X_0 \ar[ld]^{s_0} \\
& C_0 & & C_0}$$
The internal graph $d,c \colon X_1 \rightrightarrows X_0$ inherits
a structure of  groupoid from that of $\C.$ Moreover, the internal
functor $s=(s_1,s_0) \colon \X \to \C$ is a weak equivalence (it is full and faithful
by construction of $s_1,$ and it is essentially surjective because $s_0$ is a regular
epimorphism). Finally, observe that, since $X_0$ is regular projective, the
groupoid $\X$ is $\Sigma$-projective: every weak equivalence with codomain
$\X$ is in fact an equivalence.
}\end{Replacement}

\begin{Proposition}\label{PropReplacement}
Let $\C$ and $\D$ be  groupoids and fix a replacement $s \colon \X \to \C$
as in \ref{Replacement}. There is an equivalence of groupoids
$$\cB(\cC)(J(\C),J(\D)) \simeq \GrpdC(\X,\D)$$
\end{Proposition}

\begin{proof}
Since $s \colon \X \to \C$ is a weak equivalence, $\cF(s) \colon J(\X) \to J(\C)$
is an equivalence (see condition EF0 in the proof of Theorem \ref{ThMain}).
Therefore, $\cF(s)$ induces an equivalence
$$\cB(\cC)(J(\C),J(\D)) \simeq \cB(\cC)(J(\X),J(\D))$$
Moreover, since $X_0$ is regular projective, all extensions of the form
$X_0 \leftarrow E \leftarrow D$ split and then all butterflies from $J(\X)$
to $J(\D)$ split:
$$\cB(\cC)(J(\X),J(\D)) = \cB(\cC)(J(\X),J(\D))_{\rm{split}}$$
Finally, following \ref{SplitFarfalleToMor}, we have
$$\cB(\cC)(J(\X),J(\D))_{\rm{split}} \simeq \GrpdC(\X,\D)$$
\end{proof}

\begin{Remark}\label{RemGenReplacement}{\em
To end this section, we sketch a general argument on bicategories of fractions
which subsumes Proposition \ref{PropReplacement}.
\par
Let $\Sigma$ be a class of 1-cells in a bicategory $\cB$ with a right calculus of fractions. Assume that:
\begin{enumerate}
\item $\Sigma$ satisfies the $2 \Rightarrow 3$ property:
let $F \colon \C \to \D$ and $G \colon \D \to \E$ be 1-cells in $\cB;$
if two of $F,$ $G$ and $F \cdot G$ are in $\Sigma,$ then the third one is in $\Sigma;$
\item Every 1-cell $W \colon \C \to \D$ of $\Sigma$ is full and
faithful, that is, for all objects $\A$ the functor
$\cB(\A,W) \colon \cB(\A,\C) \to \cB(\A,\D)$
is full and faithful;
\item For every object $\C$ in $\cB$ there exists $S \colon \X \to \C$
in $\Sigma$ with $\X$ a $\Sigma$-projective object.
\end{enumerate}
Then, if we fix objects $\C$ and $\D$ and a 1-cell $S \colon \X \to \C$
as in 3, the functor assigning to a 1-cell $F \colon \X \to \D$ the span
$$\xymatrix{ & \X \ar[ld]_S \ar[rd]^F \\
\C & & \D}$$
yields an equivalence of groupoids
$$\cB(\X,\D) \simeq \cB[\Sigma^{-1}](\C,\D)$$
}\end{Remark}

\section{A reminder on semi-abelian categories}\label{sec:semiab}

The general context where the theory of crossed modules and weak maps can be developed
along the lines described in this paper is that of semi-abelian categories where an
additional assumption is made: the commutativity (in the sense of Smith)  of internal equivalence relations
is determined by the commutativity (in the sense of Huq) of their {\em zero}-classes. In the following, the basic notions are
recalled and the notations are fixed, for the reader's convenience.
\begin{SemiAbCat}\label{SemiAbCat}{\em
Semi-abelian categories where introduced in 2002 \cite{JMT02}, and they represent the
{\em state-of-the-art}\/ in the long-lasting investigations whose aim is to provide an
abstract categorical setting for non (necessarily) commutative pointed algebraic structures, such
as groups, rings or Lie algebras.

A category is semi-abelian when it is pointed (i.e.\ $0=1$) with finite
coproducts, protomodular \cite{Bourn} and exact (in the sense of Barr).

Pointed protomodular categories can be characterized as pointed, finitely complete categories where the
{\em split short five lemma} holds: given the diagram
$$
\xymatrix{
K\ar[r]^k \ar[d]_{h}
&A\ar[d]^f \ar@<.5ex>[r]^{p}
&B\ar@<.5ex>[l]^{s}\ar[d]^g
\\
K'\ar[r]_{k'} &A'\ar@<.5ex>[r]^{p'} &B'\ar@<.5ex>[l]^{s'}
}
$$
with $kf=hk'$, $pg=fp'$ and $sf=gs'$, the morphism $f$ is an isomorphism if $k$ and $g$ are.

Recall that a category is exact (in the sense of Barr) when it is {\em regular} and
internal equivalence relations are {\em effective}, i.e.\ kernel pairs.

Finally a \emph{regular category} is a finitely complete category where  effective equivalence relations
 have coequalizers that are stable under pullbacks.\\

The protomodularity condition can be reformulated when the category $\cC$ is pointed regular. This is stated by the
following useful characterization, that we quote from \cite{BG}, IV.4.A:
\begin{nBarrKock}\label{nBarrKock}
Let $\cC$ be a regular pointed category. Then  $\cC$ is protomodular iff
 the ``normalized Barr-Kock'' property holds: in any commutative diagram
$$
\xymatrix{
K\ar[r]^k \ar[d]_{h}
&A\ar[d]^f \ar[r]^{p}\ar@{}[dr]|{(1)}
&B\ar[d]^g
\\
K'\ar[r]_{k'} &A'\ar[r]_{p'} &B'
}
$$
where $K=\ker p$, $K'=\ker p'$ and  with $p$ being regular epi,
$$
h\ \mathrm{is\ an\ iso\ iff}\ \mathrm{the\ square\ }(1)\ \mathrm{is\ a\ pullback.}
$$
\end{nBarrKock}
}\end{SemiAbCat}
\begin{HuqIsSmith}\label{HuqIsSmith} {\em
In order to introduce the so-called ``Huq = Smith'' condition, we first recall the notions of commuting subobjects
and commuting equivalence relations.

Two subobjects $\xymatrix{G\ar[r]^g& E&\ar[l]_h H}$ commute in the sense of Huq (see \cite{Huq68,BG02})
if they cooperate as morphisms, i.e.\ if there exists a (unique) morphism $\varphi$ such that the diagram
$$
\xymatrix@C=10ex{G\ar[r]^(.4){\langle1,0\rangle}\ar[dr]_g&G\times H\ar@{-->}[d]^{\varphi}
&\ar[dl]^h H\ar[l]_(.4){\langle0,1\rangle}
\\&E}
$$
commutes.

Suppose that the maps $g,h$ are normal monomorphisms, i.e.\ kernels. Then the denormalized version of the above notion is that
of commuting equivalence relations.
A pair of equivalence relations on a common object $E$
$$
\xymatrix@C=10ex{
R \ar@<1ex>[r]^{r_1} \ar@<-1ex>[r]_{r_2} &
E \ar[l]|{e_R}\ar[r]|{e_S}&
S\ar@<1ex>[l]^{s_2} \ar@<-1ex>[l]_{s_1}
}
$$
commutes (in the sense of Smith, see \cite{Smi76, Ped95}) when there exists a (unique) morphism $\Phi$ such that
the diagram
$$
\xymatrix@C=10ex{R\ar[r]^(.4){\langle 1,r_1 e_S\rangle}\ar[dr]_{r_0}&R\times_{r_1,s_0} S\ar@{-->}[d]^{\Phi}
&\ar[dl]^{s_1} S\ar[l]_(.4){\langle s_0e_R ,1\rangle}
\\&E}
$$
commutes.
It is a well known fact that when two equivalence relations commute, so do their normalizations (see \cite{BG02}).

The converse does not hold in general, not even in semi-abelian categories (see \cite{Bourn04}
 for a counterexample, due to G. Janelidze, in the semi-abelian category of digroups). Nevertheless it does hold in several important algebraic
contexts, as pointed {\em strongly protomodular} categories (see \cite{BG02}, section 6) and pointed {\em action
accessible} categories (see \cite{MM10}). As a matter of fact, for internal structures in (many) pointed algebraic varieties,
this is quite a crucial notion and it recaptures the feeling that a local behavior near the identity element determines
a global behavior. Furthermore, it has been acknowledged in \cite{Janelidze10} that this property is a candidate
to become an axiom for ``good'' semi-abelian categories.
In present work we will refer to it as to the ``Huq = Smith'' property.
\begin{Remark}\label{rem:HuqIsSmith}
{\em
Butterflies were originally defined by B. Noohi for crossed modules of groups \cite{Noohi05}, but the author
himself, in \cite{AN09} with E. Aldrovandi, extends the construction to crossed modules
of internal groups in the Grothendieck topos  $\widehat{\textsf{S}}$, i.e.\ the topos of the sheaves over
a site $(\textsf{S},J)$ with subcanonical topology $J$.

As a matter of fact, the present setting generalizes the one  of \cite{AN09}.
 In fact more is true: our results apply to any pointed strongly protomodular
algebraic theory in a Grothendieck topos. To see this it is necessary to recollect some results from the literature.
First, in \cite{BB}, Example 4.6.3 shows that if $\mathbb{T}$ is a pointed protomodular algebraic theory, and $\mathcal{C}$ is a
regular (exact) category, then the category $\mathrm{Alg}_{\mathbb{T}}(\mathcal{C})$ of the models of $\mathbb{T}$ in $\mathcal{C}$ is
homological (exact homological).
Hence, if  $\mathcal{C}$ is exact, the missing condition for $\mathrm{Alg}_{\mathbb{T}}(\mathcal{C})$ to be
semi-abelian is its  finite cocompleteness.

Indeed  the category of models of an algebraic theory in an elementary topos $\mathcal{E}$ is
finitely cocomplete if $(i)$ the topos has a Natural Number Object, and $(ii)$ the theory is finitely
presented. Back to the situation considered here, for a Grothendieck topos $\mathcal{E}$, condition $(i)$ is
free, and condition $(ii)$ can be dropped (see \cite{BB} again, the discussion after the cited example), so that
$\mathrm{Alg}_{\mathbb{T}}(\mathcal{E})$ is semi-abelian.

Concerning the condition ``Huq = Smith'', strongly protomodular semi-abelian (i.e.\ strongly semi-abelian)
categories have this nice property, and we know from \cite{Borceux} that for a strongly protomodular
(not necessarily pointed) theory $\mathbb{T}$, and a finitely complete category $\mathcal{E}$, the category
of models $\mathrm{Alg}_{\mathbb{T}}(\mathcal{E})$ is still strongly protomodular.
This is clearly the case for a Grothendieck topos $\mathcal{E}$.

In conclusion we can state that not only our constructions and results apply to the situation described in
\cite{AN09}, but also in the context of internal Lie algebras, internal rings and other
strongly semi-abelian theories defined internally in a Grothendieck topos $\mathcal{E}$.
}
\end{Remark}
\begin{Remark}{\em
Two morphisms cooperate if their images do, and this happens precisely when their commutator is trivial,
for a suitable notion of commutator. Unfortunately, describing the many aspects of the commutator theory
involved would take us far beyond  our purposes.
The interested reader may refer to \cite{MM10}, and the bibliography therein.
}\end{Remark}
}\end{HuqIsSmith}
\begin{InternalAct}\label{InternalAct}{\em
Diverse notions of actions exist in many algebraic contexts. Most of them share the disadvantage of not being defined
intrinsically, but as set-theoretical maps satisfying certain properties. From an algebraic-categorical point of
view this is not convenient, since  those maps are difficult to deal with. This issue has been fixed
by the notion of \emph{internal action} \cite{BJ,BJK}, that expresses its full classifying power in the broad context
of semi-abelian categories.

Let $\cC$ be a finitely complete pointed category with coproducts. Then for any object $B$ in $\cC$
one can define a functor ``ker'' from the category of split epimorphisms (points) over $B$ into $\cC$
$$
\ker : Pt_{B}(\cC) \to \cC,\quad \raisebox{4ex}{\xymatrix{A\ar@<2pt>[d]^{p}\\B\ar@<2pt>[u]^{s}}}\mapsto \ker(p).
$$
This has a left adjoint:
$$
\mathrm{B+(-)} : \cC \to Pt_{B}(\cC), \quad X \mapsto \raisebox{4ex}{\xymatrix{B+X\ar@<2pt>[d]^{[1,0]}\\B\ar@<2pt>[u]^{i_B}}},
$$
The monad corresponding to this adjunction is denoted by $B\flat(-)\colon \cC\to\cC$, and, for any object A of $\cC$,
we obtain a kernel diagram:
$$
\xymatrix{B\flat A\ar[r]^{j_{B,A}}&B+A\ar[r]^(.6){[1,0]}&B.}
$$
The $B\flat(-)$-algebras are called internal $B$-actions in $\cC$.

Let us observe that in the case of groups, the object $B\flat A$ is the group generated by the \emph{formal conjugates} of elements of $A$ by elements
of $B$, i.e. by the triples of the kind $(b,a,b^{-1})$ with $b\in B$ and $a\in A$.

For any object $A$ of $\cC$, one can define a canonical conjugation action of $A$ on $A$ itself given by the composition:
$$
\chi_A:\xymatrix{A\flat A\ar[r]^{j_{A,A}}&A+A\ar[r]^(.6){[1,1]}&A}.
$$
In the category of groups, the morphism $\chi_A$ is the internal action associated to the usual conjugation in $A$:
 the realization morphism $[1,1]$ of above makes the formal conjugates of $A\flat A$ computed effectively in $A$.
Finlly observe that conjugation actions are  components of a natural transformation
$\chi\colon (-)\flat(-) \Rightarrow \mathrm{Id}_{\cC}$.
}\end{InternalAct}

\end{document}